\documentclass[a4paper,11pt]{amsart}
\usepackage{amssymb, amsmath,latexsym,amsfonts,amsbsy, amsthm,mathtools,graphicx,CJKutf8,CJKnumb,CJKulem,color,geometry}
\usepackage{mathrsfs}
\usepackage{verbatim}
\usepackage{geometry}
\usepackage{ulem}
\usepackage{indentfirst}
\usepackage{txfonts}
\def\l{\left}
\def\r{\right}
\def\f{\frac}

\def\ez{\epsilon}
\def\dz{\delta}

\def\et{\eta}
\def\lz{\lambda}
\def\Lz{\Lambda}
\def\az{\alpha}
\def\bz{\beta}
\def\sz{\sigma}
\def\gz{\gamma}

\def\rz{\rho}

\author{Lin Tang and Qian Zhang }
\title{\bf{Global $W^{2,p}$ regularity on the linearized Monge-Amp$\grave{e}$re equation with $\mathrm{VMO}$ type coefficients}\footnotetext{ \hspace{-0.65
cm} 2010 Mathematics Subject  Classification: 35J55.\\
The  research was supported  by the NNSF  (11271024)
and (11571289)  of China.}}
\date{}

\theoremstyle{plain}
\theoremstyle{plain}\newtheorem{thm}{Theorem}[section]
\theoremstyle{plain}\newtheorem{prop}{Proposition}[section]
\theoremstyle{plain}\newtheorem{cor}{Corollary}[section]
\theoremstyle{plain}\newtheorem{lem}{Lemma}[section]
\theoremstyle{plain}
\newtheorem*{thm1}{Theorem 1}
\newtheorem*{thm2}{Theorem 2}
\newtheorem*{thm3}{Theorem 3}
\numberwithin{equation}{section}

\begin{document}
\maketitle
\noindent {\bf{Abstract}.}\quad
In this paper, we establish global $W^{2,p}$ estimates for solutions of the linearized Monge-Amp$\grave{e}$re equation
$$\mathcal{L}_{\phi}u:=\mathrm{tr}[\Phi D^2 u]=f,$$
where the density of the Monge-Amp$\grave{e}$re measure $g:=\mathrm{det}D^2\phi$ satisfies a $\mathrm{VMO}$-type condition, and $\Phi:=(\mathrm{det}D^2\phi)(D^2\phi)^{-1}$ is the cofactor matrix of $D^2\phi$.

\bigskip

\section{Introduction}\label{s1}

This paper is concerned with global $W^{2,p}$ estimates for solutions of the linearized Monge-Amp$\grave{e}$re equation
\begin{equation}\label{s1:1}
\mathcal{L}_{\phi}u:=\mathrm{tr}[\Phi D^2 u]=f\quad\mathrm{in}\;\Omega,
\end{equation}
where $\phi$ is a solution of 
\begin{equation}\label{s1:2}
\mathrm{det}D^2\phi=g,\quad\quad 0<\lz\le g\le\Lz\quad\mathrm{in}\;\Omega
\end{equation}
for some constants $0<\lz\le\Lz<\infty$.

When the density $g=\mathrm{det}D^2\phi$ is continuous, the regularity of the Monge-Amp$\grave{e}$re equation \eqref{s1:2} and the linearized Monge-Amp$\grave{e}$re equation \eqref{s1:1} has been well studied. Caffarelli proved in \cite{C} that if $g$ is continuous (resp. $C^{0,\az}$), then an interior $W^{2,p}$ estimate for any $p>1$ (resp. $C^{2,\az}$ interior estimate) holds for strictly convex solutions of \eqref{s1:2}. Corresponding global $W^{2,p}$ and $C^{2,\az}$ regularity under further assumptions on $\Omega$ and the boundary data were established by Savin \cite{S2,S3} and Trudinger and Wang \cite{TW}. For the linearized equation \eqref{s1:1}, Guti$\acute{e}$rrez and Tournier derived interior $W^{2,\dz}$ estimates for small $\dz$ in \cite{GT}, and interior $W^{2,p}$ regularity for general $1<p<q$, $f\in L^q (q>n)$ was obtained in \cite{GN2}. By using Savin's localization theorem, Nam Le and Nguyen derived in \cite{NT} corresponding global $W^{2,p}$ estimates for the Dirichlet problem 
\begin{equation}\label{u}
\left\{
\begin{array}{rcl}
\mathcal{L}_{\phi}u=f&&{\mathrm{in}\;\Omega},\\
u=\varphi&&{\mathrm{on}\;\partial\Omega}.
\end{array}\right.
\end{equation}

On the other hand, in the case that $g$ is discontinuous, Huang \cite{H} proved interior $W^{2,p}$ estimates for solutions of \eqref{s1:2} when $g$ belongs to a $\mathrm{VMO}$-type space $\mathrm{VMO}_\mathrm{loc}(\Omega,\phi)$.

In this paper, we first study the interior $W^{2,p}$ estimates for solutions of \eqref{s1:1} when $g$ belongs to the $\mathrm{VMO}$-type space defined in \cite{H}. Next, we define a global $\mathrm{VMO}$-type space $\mathrm{VMO}(\Omega,\phi)$ and establish global $W^{2,p}$ estimates for solutions of \eqref{s1:1} when $g$ belongs to $\mathrm{VMO}(\Omega,\phi)$. Hence, in the interior case, we generalize the result in \cite{GN2} from the case that $g$ is continuous to the case that $g$ belongs to $\mathrm{VMO}_{\mathrm{loc}}(\Omega,\phi)$. And in the global case, we extend \cite{NT} from the case that $g$ is continuous to the case that $g$ belongs to $\mathrm{VMO}(\Omega,\phi)$.

Our first main theorem concerns interior $W^{2,p}$ estimates for solutions of \eqref{s1:1} and is stated as follows.

\begin{thm1}\label{s2:thm GN21.1}
Let $B_{\az_n}\subset\Omega\subset B_1$ be a normalized convex domain and $\phi\in C(\overline{\Omega})$ be a solution of \eqref{s1:2} with $\phi=0$ on $\partial\Omega$, 
where $g\in\mathrm{VMO}_{\mathrm{loc}}(\Omega,\phi)$.
Let $p>1, \max\{n,p\}<q<\infty$. Let $u\in C^2(\Omega)$ be a solution of \eqref{s1:1}. Let $\Omega'\Subset\Omega$, then there exists $C>0$ depending only on $n,\lz,\Lz, p, q,\Omega'$ and $g$ such that
\begin{equation*}
\|D^2 u\|_{L^p(\Omega')}\leq C\l(\|u\|_{L^\infty(\Omega)}+\|f\|_{L^q(\Omega)}\r).
\end{equation*}
\end{thm1}

The next theorem gives global $W^{2,p}$ estimates for solutions of \eqref{s1:1} when $g$ satisfies a $\mathrm{VMO}$-type small oscillation.

\begin{thm2}\label{s3:thm NT1.1}
Assume $\Omega$ is a bounded, uniform convex domain with $\partial\Omega\in C^3$. Let $\phi\in C^{0,1}(\overline{\Omega})\cap C^2(\Omega)$ be a solution of \eqref{s1:2} with $\phi=0$ on $\partial\Omega$. Let $u\in C^2(\Omega)$ be a solution of \eqref{u} with $\varphi=0$, where $f\in L^q(\Omega)$ with $n<q<\infty$. Then, for any $1<p<q$, there exist $0<\ez<1$ and $C>0$ depending only on $n,\lz,\Lz,\Omega,p,q$ such that if
$$\sup_{x\in\overline{\Omega},\,t>0}\mathrm{mosc}_{S_{\!\phi}(x,t)}\,g\leq\ez,$$
then 
$$\|u\|_{W^{2,p}(\Omega)}\leq C\|f\|_{L^q(\Omega)}.$$
\end{thm2}

Theorem \ref{s3:thm NT1.1} implies the following result concerning global $W^{2,p}$ estimates for solutions of \eqref{s1:1} when $g$ belongs to a global $\mathrm{VMO}$-type space.

\begin{thm3}\label{s3:thm NT1.2}
Assume $\Omega$ is a bounded, uniform convex domain with $\partial\Omega\in C^3$. Let $\phi\in C^{0,1}(\overline{\Omega})\cap C^2(\Omega)$ be a solution of \eqref{s1:2} with $\phi=0$ on $\partial\Omega$, where $g\in\mathrm{VMO}(\Omega,\phi)$. Let $u\in C^2(\Omega)\cap C(\overline{\Omega})$ be a solution of \eqref{u}, where $f\in L^q(\Omega)$, $\varphi\in C^2(\Omega)\cap W^{2,s}(\Omega)$ with $n<q<s<\infty$. Then, for any $1<p<q$, there exists $C>0$ depending only on $n,\lz,\Lz,\Omega,p,q,s,g$ such that
$$\|u\|_{W^{2,p}(\Omega)}\leq C\l(\|\varphi\|_{W^{2,s}(\Omega)}+\|f\|_{L^q(\Omega)}\r).$$
\end{thm3}

The spaces $\mathrm{VMO}_{\mathrm{loc}}(\Omega,\phi)$ and $\mathrm{VMO}(\Omega,\phi)$ above are defined in Section 2.

Now we indicate the main steps in the proof of the above theorems.

In the interior case, the main ingredient is the stability of the cofactor matrix of $D^2\phi$ under a $\mathrm{VMO}$-type condition of $g=\mathrm{det}D^2\phi$. For this, we use the interior $W^{2,p}$ estimates for $\phi$ in \cite{H}. Once we have this stability result, we can use similar arguments as in the case that $g$ is continuous in \cite{GN2} to prove power decay estimates for the distribution function of $D^2u$, using the power decay estimates of the distribution function of $D^2\phi$ under a $\mathrm{VMO}$-type condition of $g$ in \cite{H} instead.  

For the global case, correspondingly we first establish the stability of the cofactor matrices at the boundary. However, due to the $\mathrm{VMO}$-type property of $g$, we need to modify the boundary data of the good Monge-Amp$\grave{e}$re equation so that the maximal principle can be applied. Next we give global power decay of the distribution function of $D^2\phi$, using the interior case in \cite{H} and the covering lemma in \cite{NT}. The main difficulty in proving Theorem 2 is a localized power decay estimate for the distribution function of $D^2\phi$ at the boundary. Since $g$ is in $\mathrm{VMO}$-type spaces, there is a difficulty in using convex envelopes of functions to help estimate the density of good sets where $D^2\phi$ is bounded. Instead, we use the idea in corresponding interior case \cite{H}. Namely, we compare $\phi$ with the solution $w$ of the good Monge-Amp$\grave{e}$re equation. We apply the localized
covering lemma at the boundary and cover a neighborhood of a point at the boundary by sections of $w$. Then we use the
one-sided small power decay estimate in \cite{GN2} in each of these sections to the linearized
operator $\mathcal{L}_w$ to estimate the set where $\phi$ can be touched from above by a quadratic polynomial. With this localized power decay estimate and arguing as in \cite{NT}, it is straightforward to prove the power decay estimate for the distribution function of $D^2 u$.

The paper is organised as follows. Sections $2$ and $3$ are devoted to the interior $W^{2,p}$ estimates for solutions $u$ of \eqref{s1:1} when $g=\mathrm{det}D^2\phi$ satisfies a $\mathrm{VMO}$-type condition. In Section $2$, we first introduce some notation and give the $W^{2,p}$ estimates for $\phi$ in \cite{H}. Using these results, we then establish the stability of cofactor matrix of $D^2\phi$ under a $\mathrm{VMO}$-type condition of $g=\mathrm{det}D^2\phi$. In Section $3$, we give the power decay estimates for the distribution functions of $D^2u$. Then we give the proof of Theorem 1. Sections $4$-$6$ are devoted to the global $W^{2,p}$ estimates for solutions $u$ of \eqref{s1:1}. In Section $4$, we prove a stability result at the boundary for the cofactor matrices of $D^2\phi$ and the global power decay estimates for the distribution functions of $D^2\phi$. In Section $5$, we prove a localized power decay estimates for the distribution function of $D^2\phi$ at the boundary. Then we use this result to establish the power decay estimates for the distribution function of $D^2u$. In the last section, we give the proof of Theorems 2 and 3.

\section{Interior estimates-Stability of cofactor matrices}\label{s2}

In this section, we first introduce some notation and give the interior $W^{2,p}$ estimates of $\phi$ in \cite{H}. Then we establish the stability of cofactor matrices of $D^2\phi$ under a $\mathrm{VMO}$-type condition of $g=\mathrm{det}D^2\phi$. 

Let $\Omega$ be a convex domain in $\mathbb{R}^n$ and $\phi\in C(\overline{\Omega})$ is a strictly convex function. A section of $\phi$ centered at $x_0\in\overline{\Omega}$ with height $h$ is defined by
$$S_{\!\phi}(x_0,h):=\{x\in\overline{\Omega}:\phi(x)<\phi(x_0)+\nabla\phi(x_0)\cdot(x-x_0)+h\},$$
when $x_0\in\partial\Omega$, the term $\nabla\phi(x_0)$ is understood in the sense that
$$x_{n+1}=\phi(x_0)+\nabla\phi(x_0)\cdot(x-x_0)$$
is a supporting hyperplane for the graph of $\phi$ but for $\epsilon>0$,
$$x_{n+1}=\phi(x_0)+(\nabla\phi(x_0)+\epsilon\nu_{x_0})\cdot(x-x_0)$$
is not a supporting hyperplane, where $\nu_{x_0}$ denotes the interior unit normal to $\partial\Omega$ at $x_0$.

For any function $u$ which is differentiable at some point $x_0\in\Omega$, we always denote for simplicity
$$l_{u,x_0}:=u(x_0)+\nabla u(x_0)\cdot(x-x_0).$$

Let $\phi\in C(\overline{\Omega})$ be a solution of \eqref{s1:2} with $\phi=0$ on $\partial\Omega$. The space $\mathrm{VMO}_\mathrm{loc}(\Omega,\phi)$ is defined in \cite{H} as follows. Given a function $g\in L^1(\Omega)$, we say that $g\in \mathrm{VMO}_\mathrm{loc}(\Omega,\phi)$ if for any $\Omega'\Subset\Omega$,
$$Q_g(r,\Omega'):=\sup_{\substack{x_0\in\Omega',\\\mathrm{diam}(S_{\!\phi}(x_0,h))\leq r}}\mathrm{mosc}_{S_{\!\phi}(x_0,h)}\,g\rightarrow 0\quad\quad\mathrm{as}\;r\rightarrow 0.$$
Correspondingly, we say that $g\in\mathrm{VMO}(\Omega,\phi)$ if
$$\et_g(r,\Omega):=\sup_{\substack{x\in\overline{\Omega},\,t>0,\\\mathrm{diam}(S_{\!\phi}(x,t))\leq r}}\mathrm{mosc}_{S_{\!\phi}(x,t)}g\rightarrow 0,\quad\quad \mathrm{as}\;r\rightarrow 0.$$
Here the mean oscillation of $g$ over a measurable subset $A\subset\overline{\Omega}$ is defined by
$$\mathrm{mosc}_A g:=\fint_A |g(x)-g_A|dx,$$
where $g_A=\fint_A g dx$ denotes the average of $g$ over $A$.

There are two simple facts about the spaces defined above.

\begin{prop}\label{s2:prop VMO}
For any function $g^1$ such that $g:=(g^1)^n\in L^1(\Omega)$, the following hold:
\begin{enumerate}
\item[\rm(i)]
\begin{eqnarray*}
\left(\int_\Omega |g^1-(g^1)_\Omega|^n dx\right)^{\frac{1}{n}}\leq 2\left(\int_\Omega |g-g_\Omega| dx\right)^{\frac{1}{n}}.
\end{eqnarray*}
\item[\rm(ii)]
For any measurable subset $A\subset\Omega$, we have
\begin{eqnarray*}
\left(\int_\Omega |g^1-(g^1)_A|^n dx\right)^{\frac{1}{n}}\leq\left\{1+\left(\frac{|\Omega|}{|A|}\right)^{\frac{1}{n}}\right\}\left(\int_\Omega |g^1-(g^1)_\Omega|^n dx\right)^{\frac{1}{n}}.
\end{eqnarray*}
\end{enumerate}
\end{prop}

In this section, we always use the following assumption:\\
\noindent$\mathbf{(H)}$ $B_{\az_n}\subset\Omega\subset B_1$ is a normalized convex domain and $\phi\in C(\overline{\Omega})$ is a solution of \eqref{s1:2} with $\phi=0$ on $\partial\Omega$.

Under this assumption, we often take $w$ to be the solution of 
\begin{equation}\label{w}
\left\{
\begin{array}{rcl}
\mathrm{det}\,D^2 w=1&&{\mathrm{in}\;\Omega},\\
w=0&&{\mathrm{on}\;\partial\Omega}.
\end{array}\right.
\end{equation}

For $0<\az<1$, we define
\begin{equation*}\label{s2:GN24.41}
\Omega_\az:=\{x\in\Omega:\phi(x)<(1-\az)\min_\Omega\phi\}.
\end{equation*}

The following Lemmas \ref{s2:lem H2.1}, \ref{s2:lem H3.1}, \ref{s2:lem H6.3} and Theorem \ref{s2:thm HA(i)} were proved in \cite{H}.

\begin{lem}(See \cite[Lemma 2.1]{H}.)\label{s2:lem H2.1}
Assume that condition $\mathbf{(H)}$ holds. Then for any $\Omega'\Subset\Omega$, there exist positive constants $h_0, C$ and $q$ such that for $x_0\in\Omega'$, and $0<h\leq h_0$,
\begin{equation*}
B_{C^{-1}h}(x_0)\subset S_{\!\phi}(x_0,h)\subset B_{Ch^q}(x_0),
\end{equation*}
where $q=q(n,\lambda,\Lambda)$ and $h_0, C$ depend only on $n,\lambda,\Lambda$ and $\mathrm{dist}(\Omega',\partial\Omega)$.
\end{lem}

The maximum principle below is used to compare solutions $\phi$ of \eqref{s1:2} and $w$ of \eqref{w}, where $g=\mathrm{det}\;D^2\phi\in\mathrm{VMO}_\mathrm{loc}(\Omega,\phi)$ defined above.

\begin{lem}(See \cite[Lemma 3.1]{H}.)\label{s2:lem H3.1}
Let $\phi$ and $w$ be the weak solutions to $\mathrm{det}D^2\phi=g_1^n\geq 0$ and $\mathrm{det}D^2 w=g_2^n\geq 0$ in $\Omega$, respectively. Assume that $g_1,g_2\in L^n(\Omega)$. Then
$$\max_{\overline{\Omega}}(\phi-w)\leq \max_{\partial\Omega}(\phi-w)+C_n\mathrm{diam}(\Omega)\left(\int_\Omega(g_2-g_1)^{+n}dx\right)^{\f{1}{n}}.$$
\end{lem}

Next we give the $W^{2,p}$ estimates of solutions $\phi$ of \eqref{s1:2} in \cite{H}.

\begin{thm}(See \cite[Theorem A(i)]{H}.)\label{s2:thm HA(i)}
Assume that condition $\mathbf{(H)}$ holds. 
\begin{enumerate}
\item[(i)] Let $0<\az<1, 1\leq p<\infty$ and denote $\az_0:=\f{\az+1}{2}$. There exist constants $0<\ez<1$ and $C>0$ depending only on $n,\lz,\Lz, p$ and $\az$ such that if $\mathrm{mosc}_{S}g\leq\ez$ for any $S=S_{\!\phi}(x_0,h)\subset\Omega_{\f{\az_0+1}{2}}$, then
    $$\|D^2\phi\|_{L^p(\Omega_{\az})}\leq C.$$
\item[(ii)]If $g\in \mathrm{VMO}_\mathrm{loc}(\Omega,\phi)$, then
$D^2\phi\in L^p_{\mathrm{loc}}(\Omega)$ for any $1\leq p<\infty$.
\end{enumerate}
\end{thm}

Using Lemma \ref{s2:thm HA(i)} we are ready to prove stability of cofactor matrices of $D^2\phi$ under a $\mathrm{VMO}$-type small oscillation of $g$.

\begin{lem}\label{s2:lem GN13.4}
Let $B_{\frac{6}{5}}\subset\Omega^k\subset B_n$ be a sequence of normalized convex domain converging in the Hausdorff metric to a normalized convex domain $B_{\frac{6}{5}}\subset\Omega\subset B_n$. For each $k\in\mathbb{N}$, let $\phi_k\in C(\overline{\Omega^k})$ be a convex function satisfying
$$\left\{
\begin{array}{rcl}
\mathrm{det}D^2\phi_k=g_k&&{\mathrm{in}\;\Omega^k},\\
\phi_k=0&&{\mathrm{on}\;\partial\Omega^k},
\end{array}\right.$$
where $0<\lambda\leq g_k=(g^1_k)^n\leq\Lambda$ in $\Omega^k$, $\mathrm{mosc}_{\Omega^k}\,g_k\leq \frac{1}{k}$ and
\begin{equation}
\sup_{S_{\!\phi_k}(x,h)\Subset\Omega^k}\,\mathrm{mosc}_{S_{\!\phi_k}(x,h)}g_k\leq\f{1}{k}.
\end{equation}
Suppose that $\phi_k$ converges uniformly on compact subsets of $\Omega$ to a convex function $\phi\in C(\overline{\Omega})$ which is a solution of
$$\left\{
\begin{array}{rcl}
\mathrm{det}D^2\phi=1&&{\mathrm{in}\;\Omega},\\
\phi=0&&{\mathrm{on}\;\partial\Omega}.
\end{array}\right.$$
Then there exists a subsequence which we still denote by $\phi_{k}$ such that for any $1\leq p<\infty$,
\begin{equation*}
\lim_{k\rightarrow\infty}\|D^2\phi_{k}-(g^1_{k})_{B_1}D^2\phi\|_{L^p(B_1)}=0,
\end{equation*}
and
\begin{equation*}
\lim_{k\rightarrow\infty}\|\Phi_{k}-(g^1_{k})_{B_1}^{n-1}\Phi\|_{L^p(B_1)}=0,
\end{equation*}
where $\Phi_{k}$ and $\Phi$ are the cofactor matrices of $D^2\phi_{k}$ and $D^2\phi$ respectively.

\begin{proof}
First we note that since $\mathrm{dist}(B_1,\partial\Omega^k)\geq \mathrm{dist}(B_1,B_{\frac{6}{5}})=\frac{1}{5}$, then $B_1\subset\Omega^k_{\az}$ for some $\az$ depending only on $n,\lz,\Lz$. For any $1\leq p<\infty$, let $\ez(p)=\ez(n,\lz,\Lz,p,\az)=\ez(n,\lz,\Lz,p)$ be the constant in Theorem \ref{s2:thm HA(i)} $(i)$, then for any $k\geq k_{\ez(p)}:=\l[\f{1}{\ez(p)}\r]+1$ we have
\begin{equation}\label{s2:GN13.4*}
\sup_{S_{\!\phi_k}(x,h)\Subset\Omega^k}\,\mathrm{mosc}_{S_{\!\phi_k}(x,h)}g_k\leq\ez(p).
\end{equation}
Thus Theorem \ref{s2:thm HA(i)} $(i)$ implies that
\begin{equation}\label{s2:GN13.41}
\|D^2\phi_k\|_{L^p(B_1)}\leq\|D^2\phi_k\|_{L^p(\Omega_{\az})}\leq C(n,\lambda,\Lambda,p)\quad\forall k\geq k_{\ez(p)}.
\end{equation}

Let $\dz>0$ be an arbitrary small constant, and let $\Omega(\dz):=\{x\in\Omega:\mathrm{dist}(x,\partial\Omega)>\dz\}$. Then there exists $k_\dz\in\mathbb{N}$ such that for all $k\geq k_\dz$,
$$\mathrm{dist}(x,\partial\Omega^k)\leq 2\dz,\quad\quad\forall x\in\partial(\Omega(\dz)).$$
Then Aleksandrov's estimate (\cite[Theorem 1.4.2]{G}) implies that 
$$|\phi_k(x)-(g^1_k)_{B_1}\phi(x)|\leq C(n,\lambda,\Lambda)\dz^{\f{1}{n}}\quad\quad\forall x\in\partial(\Omega(\dz)),\forall k\ge k_\dz.$$
By choosing $k_\dz$ even larger, we have $\Omega(\dz)\subset\Omega^k$ for $k\ge k_\dz$. Then by Proposition \ref{s2:prop VMO},
\begin{eqnarray}\label{s2:GN13.45}
\left(\int_{\Omega(\dz)} |g^1_k-(g^1_k)_{B_1}|^n dx\right)^{\frac{1}{n}}&\leq&\left(\int_{\Omega^k} |g^1_k-(g^1_k)_{B_1}|^n dx\right)^{\frac{1}{n}}\leq C(n)\left(\int_{\Omega^k} |g^1_k-(g^1_k)_{\Omega^k}|^n dx\right)^{\frac{1}{n}}\nonumber\\
&\leq&C(n)\left(\int_{\Omega_k} |g_k-(g_k)_{\Omega_k}| dx\right)^{\frac{1}{n}}\leq \frac{C(n)}{k^{\f{1}{n}}}\quad\quad\forall k\ge k_\dz.
\end{eqnarray}

Using the above two estimates and applying Lemma \ref{s2:lem H3.1} with $\phi\rightsquigarrow\phi_k, w\rightsquigarrow (g^1_k)_{B_1}\phi$, we get
\begin{eqnarray}\label{s2:GN13.42}
\max_{\overline{\Omega(\dz)}}|\phi_k-(g^1_k)_{B_1}\phi|
&\leq& \max_{\partial(\Omega(\dz))}|\phi_k-(g^1_k)_{B_1}\phi|
+C_n\mathrm{diam}(\Omega(\dz))\left(\int_{\Omega(\dz)}|g^1_k-(g^1_k)_{B_1}|^{n}dx\right)^{\f{1}{n}}\nonumber\\
&\leq&C(n,\lambda,\Lambda)\dz^{\f{1}{n}}+\frac{C(n)}{k^{\f{1}{n}}} \quad\quad \forall k\ge k_\dz.
\end{eqnarray}

Using \eqref{s2:GN13.42}, \eqref{s2:GN13.41} and similar arguments to the proof of \cite[Lemma 3.4]{GN1}, we obtain the first conclusion of the lemma. For the second conclusion, we write
$$\Phi_{k}-(g^1_{k})_{B_1}^{n-1}\Phi=\l[1-\f{(g^1_{k})_{B_1}^{n}}{\mathrm{det}D^2\phi_{k}}\r]\Phi_{k}
-\f{(g^1_{k})_{B_1}^{n-1}}{\mathrm{det}D^2\phi_{k}}\Phi_{k}\l(D^2\phi_{k}-(g^1_{k})_{B_1}D^2\phi\r)\Phi.$$
For any $1\leq q,r<\infty$, if $qr\leq n$, then we have by \eqref{s2:GN13.45}
$$\l(\int_{B_1}|g^1_{k}-(g^1_{k})_{B_1}|^{qr}dx\r)^{\f{1}{qr}}\leq C(n,q,r)\l(\int_{B_1}|g^1_{k}-(g^1_{k})_{B_1}|^{n}dx\r)^{\f{1}{n}}\leq\f{C(n,q,r)}{k^{\f{1}{n}}}$$
by H$\ddot{o}$lder inequality. On the other hand, if $qr>n$ then
\begin{eqnarray*}
\l(\int_{B_1}|g^1_{k}-(g^1_{k})_{B_1}|^{qr}dx\r)^{\f{1}{qr}}&\leq& \l(\int_{B_1}|g^1_{k}-(g^1_{k})_{B_1}|^{n}|g^1_{k}-(g^1_{k})_{B_1}|^{qr-n}dx\r)^{\f{1}{qr}}\\
&\leq&C(n,q,r,\lz,\Lz)\l(\int_{B_1}|g^1_{k}-(g^1_{k})_{B_1}|^{n}dx\r)^{\f{1}{qr}}\leq\f{C(n,q,r,\lz,\Lz)}{k^{\f{1}{qr}}}.
\end{eqnarray*}
Note that
\begin{eqnarray*}
|(g^1_{k})^{n}-(g^1_{k})_{B_1}^{n}|\leq C(n,\lz,\Lz)|g^1_{k}-(g^1_{k})_{B_1}|.
\end{eqnarray*}
Therefore, 
$$\l\|1-\f{(g^1_{k})_{B_1}^{n}}{\mathrm{det}D^2\phi_{k}}\r\|_{L^{qr}(B_1)}\to 0\quad \mathrm{as}\;k\to\infty.$$
The rest of the proof is similar to \cite[Lemma 3.5]{GN1}, using \eqref{s2:GN13.41} and the first conclusion of the lemma instead.
\end{proof}
\end{lem}

Using Lemma \ref{s2:lem GN13.4} and Blaschke selection theorem, we can obtain: 

\begin{cor}\label{s2:cor GN24.4}
Given $0<\ez_0<1$ and $0<\az<1$. Then there exists $\ez>0$ depending only on $\ez_0,n,\lz,\Lz,\az$ such that if the assumption $\mathbf{(H)}$ holds and $w\in C(\overline{\Omega})$ is the solution of \eqref{w}, where 
$$\mathrm{mosc}_{\Omega}\,g\leq \ez\quad\quad\mathrm{and}\quad\quad\sup_{S_{\!\phi}(x,h)\Subset\Omega}\mathrm{mosc}_{S_{\!\phi}(x,h)}\,g\leq\ez,$$
then
\begin{equation*}
\|\Phi-(g^1)_{\Omega_\az}^{n-1}W\|_{L^n(\Omega_\az)}\leq\ez_0,
\end{equation*}
where $\Phi$ and $W$ are the cofactor matrices of $D^2\phi$ and $D^2 w$ respectively.
\end{cor}

Next we begin to prove the power decay of the distribution function of $D^2u$. The following lemma was given in \cite{GN2}.

\begin{lem}(See \cite[Lemma 2.7]{GN2}.)\label{s2:lem GN22.7}
Assume that condition $\mathbf{(H)}$ holds. Let $u\in C^2(\Omega)$. Then for $\kappa>1$, we have
\begin{equation}\label{s2:GN22.72}
\{x\in\Omega_{\az}:|D_{i j}u(x)|>\bz^\kappa\}\subset\l(\Omega_{\az}\backslash D^{\az}_{(c\bz^{\f{\kappa-1}{2}})^{\f{1}{n-1}}}\r)\cup\l(\Omega_{\az}\backslash G_{\bz}(u,\Omega)\r)
\end{equation}
for any $\bz>0$ satisfying $(c\bz^{\f{\kappa-1}{2}})^{\f{1}{n-1}}\geq\f{\mathrm{diam}\,\Omega}{\sqrt{\et_{\az}}}$, where $c>0$ is a constant depending only on $n,\lz,\Lz$.
\end{lem}

We recall the related definitions in the above lemma. 

Let $\Omega$ be a bounded convex set in $\mathbb{R}^n$ and $\phi\in C^1(\Omega)$ be a convex function. For $u\in C(\Omega)$ and $M>0$, define the sets
\begin{eqnarray*}
G^+_{\!M}(u,\Omega)=\l\{\bar{x}\in\Omega: u\;\mathrm{is\;differentiable\;at\;}\bar{x}\;\mathrm{and}\,
u(x)\leq l_{u,\bar{x}}(x)+\f{M}{2}d(x,\bar{x})^2,\;\forall x\in\Omega\r\};\\
G^-_{\!M}(u,\Omega)=\l\{\bar{x}\in\Omega: u\;\mathrm{is\;differentiable\;at\;}\bar{x}\;\mathrm{and}\,
u(x)\geq l_{u,\bar{x}}(x)-\f{M}{2}d(x,\bar{x})^2,\;\forall x\in\Omega\r\};
\end{eqnarray*}
and $G_{\!M}(u,\Omega):=G^+_{\!M}(u,\Omega)\cap G^-_{\!M}(u,\Omega)$, where for any $x\in\Omega$ and $x_0\in\Omega$,
$$d(x,x_0)^2:=\phi(x)-l_{\phi,x_0}(x).$$

Under the assumption $\mathbf{(H)}$, for $0<\az<1$, let $\et_{\az}$ be the constant such that the conclusion in Lemma \ref{s2:lem H2.1} holds for $\Omega'\rightsquigarrow\Omega_{\az}$ and $S_{\!\phi}(x,h)\subset\Omega_{(\az+1)/2}$ for $x\in\Omega_{\az}$ and $h\leq\et_{\az}$. Denote
$$D^{\az}_{\bz}=\{x\in\Omega_{\az}:S_{\!\phi}(x,h)\subset B_{\bz\sqrt{h}}(x)\;\mathrm{for}\;h\leq\et_{\az}\}$$
and
\begin{equation}\label{s2:GN22.71}
A_{\sz}=\{\bar{x}\in\Omega:\phi(x)-l_{\phi,\bar{x}}(x)\geq\sz|x-\bar{x}|^2,\;\forall x\in\Omega\}
\end{equation}
By \cite[Lemma 6.2.2]{G}, $D^{\az}_{\bz}=\Omega_{\az}\cap A_{\bz^{-2}}$ if $\bz\geq\f{\mathrm{diam}\,\Omega}{\sqrt{\et_{\az}}}$.

The estimate of the first term on the right-hand side of \eqref{s2:GN22.72} is given by the following lemma. It is an easy corollary of Theorem 6.3 and Equation (6.12) in \cite{H}.

\begin{lem}\label{s2:lem H6.3}
Let $0<\az<1$ and denote $\az_0:=\f{\az+1}{2}$. Assume that condition $\mathbf{(H)}$ holds, where $\mathrm{mosc}_{S} g\leq\ez$ with small $\ez$ for any section $S=S_{\!\phi}(x,h)\subset\Omega_{\f{\az_0+1}{2}}$. Then there exist constants $M>0$ depending only on $n,\lz,\Lz,\az$ and $0<\dz_1<1$ depending only on $n,\lz,\Lz$ such that
$$|\Omega_{\az}\backslash D^{\az}_{s}|\leq\f{|\Omega|}{(2\ez^{\dz_1/2})^2}s^{\f{\ln\l(2\ez^{\dz_1/2}\r)}{\ln M}}$$
for each $s>0$.
\end{lem}

\section{Interior estimates-Power decay results}\label{s3}

In this section, we give the estimate of the second term on the right-hand side of \eqref{s2:GN22.72} and the proof of Theorem 1. We first include the following result from \cite{H}. Recall the definition of $A_{\sz}$ in \eqref{s2:GN22.71}.

\begin{lem}(See \cite[Lemma 6.1]{H}.)\label{s2:lem H6.1}
Assume that condition $\mathbf{(H)}$ holds, where $\mathrm{mosc}_\Omega g\leq\ez$. Let $0<\az<1$. Then there exist $0<\dz_1<1$ and $\sz>0$ depending only on $n,\lz,\Lz$ and $\az$ such that
$$|\Omega_{\az}\backslash A_{\sz}|\leq\ez^{\dz_1}|\Omega_{\az}|.$$
\end{lem}

Next, we compare explicitly two solutions originating from two different linearized Monge-Amp$\grave{e}$re equations. The following lemma is a slight modification of \cite[Lemma 4.1]{GN2}.

\begin{lem}\label{s2:lem GN24.1}
Assume that condition $\mathbf{(H)}$ holds and $w$ is the solution of \eqref{w} with $\Omega\rightsquigarrow U$ and denote $g^1:=g^{\f{1}{n}}$. Let $u\in W^{2,n}_{\mathrm{loc}}(U)\cap C(\overline{U})$ be a solution of $\mathcal{L}_\phi u=f$ in $U$ with $|u|\leq 1$ in $U$. Assume $0<\az_1<1$ and $h\in W^{2,n}_{\mathrm{loc}}(U_{\az_1})\cap C(\overline{U_{\az_1}})$ is a solution of
$$\left\{
\begin{array}{rcl}
\mathcal{L}_wh=0&&{\mathrm{in}\;U_{\az_1}},\\
h=u&&{\mathrm{on}\;\partial U_{\az_1}}.
\end{array}\right.$$
Then there exists $0<\gz<1$ depending only on $n,\lz,\Lz$ such that for any $0<\az_2<\az_1$, we have
\begin{eqnarray*}
\|u-h\|_{L^\infty(U_{\az_2})}+\l\|f-\mathrm{tr}\l(\l[\Phi-(g^1)^{n-1}_{U_{\az_1}}W\r] D^2 h\r)\r\|_{L^n(U_{\az_2})}
\leq C\l\{\|\Phi-(g^1)^{n-1}_{U_{\az_1}}W\|^{\gz}_{L^n(U_{\az_1})}+\|f\|_{L^n(U)}\r\}
\end{eqnarray*}
provided that $\|\Phi-(g^1)^{n-1}_{U_{\az_1}}W\|_{L^n(U_{\az_1})}\leq(\az_1-\az_2)^{\f{2 n}{1+(n-1)\gz}}$, where $C=C(n,\lz,\Lz,\az_1,\az_2)$.
\end{lem}

Using Lemma \ref{s2:lem H6.1} and following similar lines as in the proof of \cite[Lemma 4.2]{GN2}, we obtain the density estimate below.

\begin{lem}\label{s2:lem GN24.2}
Let $0<\ez<1, 0<\az_0<1$. Assume $\Omega$ is a bounded convex set,  $U\subset\Omega$ and the condition $\mathbf{(H)}$ holds with $\Omega\rightsquigarrow U$, where $\phi\in C^1(\Omega)\cap W^{2,n}_{\mathrm{loc}}(U)$ and $\mathrm{mosc}_U g\leq\ez$. Denote $g^1:=g^{\f{1}{n}}$. Let $u\in C(\Omega)\cap W^{2,n}_{\mathrm{loc}}(U)\cap C^1(U)$ be a solution of $\mathcal{L}_\phi u=f$ in $U$ with $|u|\leq 1$ in $U$ and
$$|u(x)|\leq C^*[1+d(x,x_0)^2]\quad\quad\mathrm{in}\;\Omega\backslash U$$
for some $x_0\in U_{\az_0}$. Then for any $0<\az\leq\az_0$, there exist constants $C,\tau>0$ and $0<\dz_1<1$ depending only on $n,\lz,\Lz,\az_0,\az$ such that
$$|G_N(u,\Omega)\cap U_{\az}|\geq\l\{1-C\l(N^{-\tau}\dz_0^{\tau}+\ez^{\dz_1}\r)\r\}|{U_\az}|$$
for any $N\geq N_0=N_0(n,\lz,\Lz,\az_0,\az,C^*)$ and provided that $\|\Phi-(g^1)^{n-1}_{U_{\f{\az_0+1}{2}}}W\|_{L^n(U_{\f{\az_0+1}{2}})}\leq\l(\f{1-\az_0}{4}\r)^{\f{2 n}{1+(n-1)\gz}}$. Here $W,\gz$ are from Lemma \ref{s2:lem GN24.1} and
$$\dz_0:=\l(\fint_{U_{\f{\az_0+1}{2}}}\|\Phi-(g^1)^{n-1}_{U_{\f{\az_0+1}{2}}}W\|^n dx\r)^{\f{\gz}{n}}
+\l(\fint_U |f|^n dx\r)^{\f{1}{n}}.$$
\end{lem}

By using Lemma \ref{s2:lem GN24.2}, a localization process and stability of cofactor matrices (Corollary \ref{s2:cor GN24.4}), we obtain the following two lemmas.

\begin{lem}\label{s2:lem GN24.3}
Let $0<\ez_0<1, 0<\az_0<1$. Let $B_{\az_n}\subset\Omega\subset B_1$ be a normalized convex domain and $\Omega_1$ be a bounded convex set such that $\Omega\subset\Omega_1$. Assume $\phi\in C^2(\Omega_1)$ is a convex function satisfying
$$\left\{
\begin{array}{rcl}
\mathrm{det}D^2\phi=g&&{\mathrm{in}\;\Omega_1},\\
\phi=0&&{\mathrm{on}\;\partial\Omega}.
\end{array}\right.$$
where $0<\lambda\leq g=(g^1)^n\leq\Lambda$ in $\Omega_1$. Let $u\in C(\Omega_1)\cap C^1(\Omega)\cap W^{2,n}_{\mathrm{loc}}(\Omega)$ be a solution of $\mathcal{L}_\phi u=f$ in $\Omega$ with $\|u\|_{L^\infty(\Omega_1)}\leq 1$. There exists $\ez>0$ depending only on $n,\lz,\Lz,\ez_0,\az_0$ such that if $S_{\!\phi}\l(x_0,\f{t_0}{\az_0}\r)\subset\Omega_{\f{\az_0+1}{2}}$ and
$$\sup_{S_{\!\phi}(x,h)\subseteq S_{\!\phi}\l(x_0,\f{t_0}{\az_0}\r)}\,g\leq\ez,$$
then
\begin{equation}\label{s2:GN24.31}
|G_{\f{N}{t_0}}(u,\Omega_1)\cap S_{\!\phi}(x_0,t_0)|\geq\l\{1-\ez_0-C\l(\f{t_0}{N}\r)^\tau\l(\fint_{S_{\!\phi}\l(x_0,\f{t_0}{\az_0}\r)}|f|^n dx\r)^{\f{\tau}{n}}\r\}|S_{\!\phi}(x_0,t_0)|
\end{equation}
for any $N\geq N_0$. Here $C,\tau,N_0$ are positive constants depending only on $n,\lz,\Lz,\az_0$.

\begin{proof}
Let $T$ be an affine map such that
\begin{equation*}\label{s2:GN24.32}
B_{\az_n}\subset T\l(S_{\!\phi}\l(x_0,\f{t_0}{\az_0}\r)\r)\subset B_1.
\end{equation*}

Let $U:=T\l(S_{\!\phi}\l(x_0,\f{t_0}{\az_0}\r)\r)$. For each $y\in T(\Omega_1)$, define
$$\tilde{\phi}(y)=|\mathrm{det}\,T|^{\f{2}{n}}\l[\phi(T^{-1}y)-l_{\phi,x_0}(T^{-1}y) -\f{t_0}{\az_0}\r],\;\mathrm{and}\;\tilde{u}(y)=u(T^{-1}y).$$
Define $\tilde{g}(y):=g(T^{-1}y)$, and then
$$\mathrm{mosc}_{U}\,\tilde{g}=\mathrm{mosc}_{S_{\!\phi}\l(x_0,\f{t_0}{\az_0}\r)}\,g\leq\ez$$
and
$$\sup_{S_{\!\tilde{\phi}}(y,h)\Subset U}\,\tilde{g}
\leq\sup_{S_{\!\phi}(T^{-1}y,|\mathrm{det}\,T|^{\f{-2}{n}}h)\Subset S_{\!\phi}\l(x_0,\f{t_0}{\az_0}\r)}\mathrm{mosc}_{S_{\!\phi}(T^{-1}y,|\mathrm{det}\,T|^{\f{-2}{n}}h)}\,g\leq\ez.$$
We apply Lemma \ref{s2:lem GN24.2} to $\tilde{u}$. Then by Corollary \ref{s2:cor GN24.4} and arguing similarly as in \cite[Lemma 4.3]{GN2}, we obtain \eqref{s2:GN24.31}.
\end{proof}
\end{lem}

Similarly, using Lemma \ref{s2:lem GN24.2} and Corollary \ref{s2:cor GN24.4}, we can obtain the variant of \cite[Lemma 4.5]{GN2}.

\begin{lem}\label{s2:lem GN24.5}
Let $0<\ez_0<1, 0<\az_0<1$. Let $B_{\az_n}\subset\Omega\subset B_1$ be a normalized convex domain and $\Omega_1$ be a bounded convex set such that $\Omega\subset\Omega_1$. Assume $\phi\in C^2(\Omega_1)$ is a convex function satisfying
$$\left\{
\begin{array}{rcl}
\mathrm{det}D^2\phi=g&&{\mathrm{in}\;\Omega_1},\\
\phi=0&&{\mathrm{on}\;\partial\Omega}.
\end{array}\right.$$
where $0<\lambda\leq g=(g^1)^n\leq\Lambda$ in $\Omega_1$. Let $u\in C(\Omega_1)\cap C^1(\Omega)\cap W^{2,n}_{\mathrm{loc}}(\Omega)$ be a solution of $\mathcal{L}_\phi u=f$ in $\Omega$. There exists $\ez>0$ depending only on $n,\lz,\Lz,\ez_0,\az_0$  such that if $S_{\!\phi}(x_0,t_0)\subset\Omega_{\az_0}$ with $t_0\leq \et_{\az_0}$,
$$\sup_{S_{\!\phi}(x,h)\subseteq S_{\!\phi}\l(x_0,\f{t_0}{\az_0}\r)}\mathrm{mosc}_{S_{\!\phi}(x,h)}\,g\leq\ez,$$
and $S_{\!\phi}(x_0,t_0)\cap G_\gz(u,\Omega_1)\neq\emptyset$, then we have
\begin{equation*}\label{s2:GN24.51}
|G_{N\gz}(u,\Omega_1)\cap S_{\!\phi}(x_0,t_0)|\geq\l\{1-\ez_0-C(N\gz)^{-\tau}\l(\fint_{S_{\!\phi}\l(x_0,\f{t_0}{\az_0}\r)}|f|^n dx\r)^{\f{\tau}{n}}\r\}|S_{\!\phi}(x_0,t_0)|
\end{equation*}
for any $N\geq N_0$. Here $C,\tau,N_0,\et_{\az_0}$ are positive constants depending only on $n,\lz,\Lz,\az_0$.
\end{lem}

Using Lemmas \ref{s2:lem GN24.3}, \ref{s2:lem GN24.5}, \ref{s2:lem GN22.7} and \ref{s2:lem H6.3}, we can obtain

\begin{thm}\label{s2:thm GN24.6}
Given $0<\lz\leq\Lz$. Let $p>1, \max\{n,p\}<q<\infty$ and $0<\az<1$. Denote $\az_0:=\f{\az+1}{2}$. Assume the condition $\mathbf{(H)}$ holds. Let $u\in C^2(\Omega)$ be a solution of $\mathcal{L}_\phi u=f$ in $\Omega$. There exist $0<\ez<1$ and $C>0$ depending only on $n,\lz,\Lz,\az, p, q$ such that if
$$\sup_{S_{\!\phi}(x,h)\subset \Omega_{\f{\az_0+1}{2}}}\mathrm{mosc}_{S_{\!\phi}(x,h)}\,g\leq\ez,$$
then we have
\begin{equation*}
\|D^2 u\|_{L^p(\Omega_\az)}\leq C\l(\|u\|_{L^\infty(\Omega)}+\|f\|_{L^q(\Omega)}\r).
\end{equation*}
\end{thm}

Now we are ready to prove the interior $W^{2,p}$ estimates for solutions to \eqref{s1:1}.\\

\noindent$\mathbf{Proof\;of\;Theorem\;1:}$
By Lemma \ref{s2:lem H2.1}, there exist positive constants $h_0, C$ and $q$ depending only on $n,\lambda,\Lambda$ and $\mathrm{dist}(\Omega',\partial\Omega)$ such that
\begin{equation}\label{s2:GN21.11}
B_{C^{-1}2 h_0}(x_0)\subset S_{\!\phi}(x_0,2 h_0)\subset B_{C(2 h_0)^q}(x_0).
\end{equation}
Choose $h_0$ smaller and we can assume $S_{\!\phi}(x_0,2 h_0)\subset B_{C(2 h_0)^q}(x_0) \subset\Omega''\Subset\Omega$. Since $g\in \mathrm{VMO}_\mathrm{loc}(\Omega,\phi)$, we have
$$\et_g(r,\Omega''):=\sup_{\substack{S_{\!\phi}(x,h)\subset\Omega'',\\\mathrm{diam}(S_{\!\phi}(x,h))\leq r}}\mathrm{mosc}_{S_{\!\phi}(x,h)}\,g\rightarrow 0,\quad \mathrm{as}\;r\rightarrow 0.$$
Let $\epsilon>0$ be the constant in Theorem \ref{s2:thm GN24.6} corresponding to $\az=\f{1}{2}$, there exists $0<r_0<1$ such that $\et_g(r_0,\Omega'')<\ez$. Take $h_0$ smaller such that $\mathrm{diam}(B_{C(2 h_0)^q}(x_0))\leq r_0$, then for any $S_{\!\phi}(x,h)\subset S_{\!\phi}(x_0,2 h_0)$, we have $S_{\!\phi}(x,h)\subset\Omega''$ and $\mathrm{diam}(S_{\!\phi}(x,h))\leq r_0$, thus,
$$\mathrm{mosc}_{S_{\!\phi}(x,h)}\,g\leq\et_g(r_0,\Omega'')<\ez.$$

Let $T$ be an affine map such that
\begin{equation*}
B_{\az_n}\subset T\l(S_{\!\phi}\l(x_0,2 h_0\r)\r)\subset B_1.
\end{equation*}
By \eqref{s2:GN21.11} we have
\begin{equation}\label{s2:GN21.12}
\|T\|\leq C h^{-1}_0,\quad \|T^{-1}\|\leq C h_0^q\leq 1
\end{equation}
if $h_0$ is even smaller. Let $\widetilde{\Omega}:=T\l(S_{\!\phi}\l(x_0,2 h_0\r)\r)$. For each $y\in\widetilde{\Omega}$, define
$$\tilde{\phi}(y)=|\mathrm{det}\,T|^{\f{2}{n}}\l[\phi(T^{-1}y)-l_{\phi,x_0}(T^{-1}y)-2 h_0\r],\;\mathrm{and}\;\tilde{u}(y)=u(T^{-1}y).$$

We have
$$\mathrm{det}D^2\tilde{\phi}=\tilde{g}(y)=(\tilde{g}^1)^n\;\mathrm{in}\;\widetilde{\Omega},
\quad\quad\tilde{\phi}=0\;\mathrm{on}\;\partial\widetilde{\Omega},$$
$$\lz\leq\tilde{g}(y)=g(T^{-1}y)\leq\Lz\;\mathrm{in}\;\widetilde{\Omega},$$
$$\mathcal{L}_{\tilde{\phi}}\tilde{u}(y)=\tilde{f}(y):=|\mathrm{det}\,T|^{\f{-2}{n}}f(T^{-1}y)\;\mathrm{in}\;\widetilde{\Omega}.$$
Moreover, choose $\az=1/2$ and denote $\az_0:=\f{\az+1}{2}$, then we have
\begin{eqnarray*}
\sup_{S_{\!\tilde{\phi}}(y,h)\subset \widetilde{\Omega}_{\f{\az_0+1}{2}}}\mathrm{mosc}_{S_{\!\tilde{\phi}}(y,h)}\,\tilde{g}
&\leq&\sup_{S_{\!\tilde{\phi}}(y,h)\Subset \widetilde{\Omega}}\mathrm{mosc}_{S_{\!\tilde{\phi}}(y,h)}\,\tilde{g}\\
&\leq&\sup_{S_{\!\phi}(T^{-1}y,|\mathrm{det}\,T|^{\f{-2}{n}}h)\Subset S_{\!\phi}\l(x_0,2 h_0\r)}\mathrm{mosc}_{S_{\!\phi}(T^{-1}y,|\mathrm{det}\,T|^{\f{-2}{n}}h)}\,g\leq\ez.
\end{eqnarray*}
Applying Theorem \ref{s2:thm GN24.6} we have
\begin{equation*}
\|D^2\tilde{u}\|_{L^p(\widetilde{\Omega}_{1/2})}\leq C\l(\|\tilde{u}\|_{L^\infty(\widetilde{\Omega})}+\|\tilde{f}\|_{L^q(\tilde{\Omega})}\r).
\end{equation*}
Back to $u$, note \eqref{s2:GN21.12},
\begin{equation*}
\|D^2 u\|_{L^p(S_{\!\phi}(x_0,h_0))}\leq C\l(\|u\|_{L^\infty(\Omega)}+\|f\|_{L^q(\Omega)}\r).
\end{equation*}
From this, by a standard covering argument, we obtain the conclusion.

\section{Boundary estimate-A stability result at the boundary for the cofactor matrices}\label{s5}

In this section, we prove a stability result at the boundary for the cofactor matrices of $D^2\phi$ and then establish the global power decay estimate for the distribution function of $D^2\phi$. We fix constants $0<\lz\le\Lz<\infty,\rz>0$ and refer to all positive constants depending only $n,\lz,\Lz$ and $\rz$ as universal constants.

Assume
\begin{equation}\label{s3:1}
B_\rz(\rz e_n)\subset\Omega\subset \{x_n\geq 0\}\cap B_{\f{1}{\rz}}.
\end{equation}
\begin{equation}\label{s3:2}
\mathrm{\Omega\;contains\;an\;interior\;ball\;of\;radius\;\rz\;
tangent\;to\;\partial\Omega\;at\;each\;point\;on\;\partial\Omega\cap\;B_\rz}.
\end{equation}
Let $\phi:\overline{\Omega}\rightarrow\mathbb{R}$, $\phi\in C^{0,1}(\overline{\Omega})\cap C^2(\Omega)$ be a convex function satisfying
\begin{equation}\label{s3:3}
\mathrm{det}\,D^2\phi=g,\quad\quad 0<\lambda\leq g\leq\Lambda\quad\mathrm{in}\;\Omega.
\end{equation}
Assume further that on $\partial\Omega\cap B_\rz$, $\phi$ separates quadratically from its tangent planes on $\partial\Omega$, namely, for any $x_0\in\partial\Omega\cap B_\rz$ we have
\begin{equation}\label{s3:4}
\rz|x-x_0|^2\leq\phi(x)-\phi(x_0)-\nabla\phi(x_0)\cdot (x-x_0)\leq\rz^{-1}|x-x_0|^2,\quad\quad\forall x\in\partial\Omega.
\end{equation}

\begin{thm}(Localization Theorem \cite{S1,S2})\label{s3:thm NT2.2}
Assume $\Omega$ satisfies \eqref{s3:1} and $\phi$ satisfies \eqref{s3:3}, and
$$\phi(0)=\nabla\phi(0)=0,\quad\quad
\rz|x|^2\leq\phi(x)\leq\rz^{-1}|x|^2,\;\mathrm{on}\;\partial\Omega\cap\{x_n\leq\rz\}.$$
Then there exists a universal constant $k>0$ such that for each $h\leq k$, there is an ellipsoid $E_h$ of volume $|B_1|h^{n/2}$ satisfying
$$k E_h\cap\overline{\Omega}\subset S_{\!\phi}(0,h)\subset k^{-1}E_h\cap\overline{\Omega}.$$
Moreover, the ellipsoid $E_h$ is obtained from the ball of radius $h^{\f{1}{2}}$ by a linear transformation $A_h^{-1}$ (sliding along the $x_n=0$ plane)
\begin{eqnarray*}
\mathrm{det}A_h=1,\quad A_h x=x-\tau_h x_n,\quad \tau_h\cdot e_n=0,\\
h^{-\f{1}{2}}A_h E_h=B_1,\quad |\tau_h|\leq k^{-1}|\log h|.
\end{eqnarray*}
\end{thm}

\begin{prop}(see \cite[Proposition 3.2]{S1})\label{s3:prop NT2.3}
Let $\phi$ and $\Omega$ satisfy the hypotheses of the Localization Theorem \ref{s3:thm NT2.2}. Assume that for some $y\in\Omega$ the section $S_{\!\phi}(y,h)\subset\Omega$ is tangent to $\partial\Omega$ at $0$ for some $h\leq c$ with $c$ universal. Then 
$$\nabla\phi(y)=a e_n\quad\quad\mathrm{for\;some\;}a\in [k_0 h^{\f{1}{2}},k_0^{-1}h^{\f{1}{2}}],$$
$$k_0 E_h\subset S_{\!\phi}(y,h)-y\subset k_0^{-1}E_h,\quad\quad k_0 h^{\f{1}{2}}\leq\mathrm{dist}(y,\partial\Omega)\leq k_0^{-1}h^{\f{1}{2}},$$
where $E_h$ is the ellipsoid defined in the Localization Theorem \ref{s3:thm NT2.2}, and $k_0>0$ a universal constant.
\end{prop}

Under the assumptions of Theorem \ref{s3:thm NT2.2}, we have for all $h$ small,
\begin{equation}\label{s3:4'}
\overline{\Omega}\cap B_{h^{\f{2}{3}}}^+\subset\overline{\Omega}\cap B_{c_1 h^{\f{1}{2}}/|\log h|}^+\subset S_{\!\phi}(0,h)\subset B_{C_1 h^{\f{1}{2}}|\log h|}\subset B_{h^{\f{1}{3}}},
\end{equation}
and 
\begin{equation}\label{s3:4''}
c_1|x|^2|\log|x||^{-2}\leq\phi(x)\leq C_1|x|^2|\log|x||^{2}
\end{equation}
for any $x\in B_{c_1}$, where $c_1,C_1$ are universal (see Equation (4.3) in \cite{S1}).

Denote the rescaled function of $\phi$ and the rescaled domain of $\Omega$ by
\begin{eqnarray}\label{s3:5}
\phi_h(x):=\f{\phi(h^{\f{1}{2}}A_h^{-1}x)}{h}\quad\quad\mathrm{and}\quad\quad\Omega_h:=h^{-\f{1}{2}}A_h\Omega.
\end{eqnarray}
Denote $U_h:=S_{\!\phi_h}(0,1)=\{x\in\overline{\Omega_h}:\phi_h(x)<1\}$, then
$$\mathrm{det}D^2\phi_h=g_h=(g^1_h)^n,\quad\quad 0<\lambda\leq g_h(x):=g(h^{\f{1}{2}}A_h^{-1}x)\leq\Lambda\quad\mathrm{in}\;\Omega_h,$$
\begin{eqnarray}\label{s3:6}
\overline{\Omega_h}\cap B_k\subset h^{-\f{1}{2}}A_h S_{\!\phi}(0,h)=U_h\subset B_{k^{-1}}^+.
\end{eqnarray}

To establish the stability of the cofactor matrix of $D^2\phi$ at the boundary we compare $\phi_h$ and the solution $w$ of $\mathrm{det}\,D^2w=1$. Due to the $\mathrm{VMO}$-type condition, we change the boundary data of $w$ as follows:

Given a constant $A>0$ such that $c(n,\lz,\Lz)\leq A\leq C(n,\lz,\Lz)$. Let $w\in C(\overline{U_h})$ be the convex solution of
\begin{equation}\label{s3:NT2.71}
\left\{
\begin{array}{rcl}
\mathrm{det}D^2 w=1&&\mathrm{in}\;U_h,\\
w=A\cdot\phi&&\mathrm{on}\;\partial U_h.
\end{array}\right.
\end{equation}

Correspondingly, we modify the class $\mathcal{P}_{\lz,\,\Lz,\,\rz,\,\kappa,\,\az}$ in \cite{NT}. Fix $n,\rz,\lambda,\Lambda,\kappa$ and $\az$, the class $\mathcal{P}_{\lz,\,\Lz,\,\rz,\,\kappa,\,\az}$ and $\mathcal{P}_{\lz,\,\Lz,\,\rz,\,\kappa,*}$ consist of the quadruples $(\Omega,\phi, g, U)$ satisfying the following conditions $(i)$-$(vii)$ and $(i)$-$(vi)$ respectively:

\noindent(i) $0\in\partial\Omega,\,U\subset\Omega\subset\mathbb{R}^n$ are bounded convex domains such that
$$B_k^+\cap\overline{\Omega}\subset\overline{U}\subset B_{k^{-1}}^+\cap\overline{\Omega}.$$
\noindent(ii) $\phi:\overline{\Omega}\rightarrow\mathbb{R}^+$ is convex satisfying $\phi=1\;\mathrm{on}\;\partial U\cap\Omega$ and
$$\mathrm{det}D^2\phi=g,\quad\lambda\leq g=(g^1)^n\leq\Lambda\;\mathrm{in}\;\Omega,$$
$$\phi(0)=0,\quad\nabla\phi(0)=0,\quad\partial\Omega\cap\{\phi<1\}=\partial U\cap\{\phi<1\}.$$
\noindent(iii)\,(quadratic separation)
$$\f{\rz}{4}|x-x_0|^2\leq\phi(x)-\phi(x_0)-\nabla\phi(x_0)\cdot (x-x_0)\leq \f{4}{\rz}|x-x_0|^2,\quad\quad\forall x,x_0\in\partial\Omega\cap B_{\f{2}{k}}.$$
\noindent(iv)\,(flatness)
$$\partial\Omega\cap\{\phi<1\}\subset G\subset\{x_n\leq\kappa\},$$
where $G\subset B_{2/k}$ is a graph in the $e_n$ direction and its $C^{1,1}$ norm is bounded by $\kappa$.\\
\noindent(v)\,(localization and gradient estimates) $\phi$ satisfies in $U$ the hypotheses of the Localization Theorem in \ref{s3:thm NT2.2} at all points on $\partial U\cap B_{c}$ and
$$|\nabla\phi|\leq C_0\quad\mathrm{in}\;U\cap B_c.$$
\noindent(vi)\,(maximal sections around the origin) If $y\in U\cap B_{c^2}$ then the maximal interior section of $\phi$ in $U$ satisfies
$$k_0^2\mathrm{dist}^2(y,\partial U)\leq\bar{h}(y)\leq c\quad\quad\mathrm{and}\quad\quad S_\phi(y,\bar{h}(y))\subset U\cap B_c.$$
\noindent(vii)\,(Pogorelov estimates)
$$\|\partial U\cap B_c\|_{C^{2,\az}}\leq c_0^{-1}$$
and if $w$ is the convex solution to
\begin{equation}\label{s3:P}
\left\{
\begin{array}{rcl}
\mathrm{det}D^2 w=1&&{\mathrm{in}\;U},\\
w=A\cdot\phi&&{\mathrm{on}\;\partial U}.
\end{array}\right.
\end{equation}
for some constant $A>0$ such that $c(n,\lz,\Lz)\leq A\leq C(n,\lz,\Lz)$, then
$$\|w\|_{C^{2,\az}(\overline{U\cap B_c})}\leq c_0^{-1}\quad\mathrm{and}\quad c_0I_n\leq D^2 w\leq c_0^{-1}I_n\;\mathrm{in}\;U\cap B_c.$$
The constants $k, k_0, c, C_0$ above depend only on $n,\lambda,\Lambda$ and $\rz$ and $c_0$ depends also on $\az$.

The class $\mathcal{P}_{\lz,\,\Lz,\,\rz,\,\kappa,*}$ is the same as that in \cite{NT}, while $\mathcal{P}_{\lz,\,\Lz,\,\rz,\,\kappa,\,\az}$ is slightly different from that in \cite{NT}.

Let $\Omega,\phi$ and $g$ satisfy \eqref{s3:1}-\eqref{s3:4}. It follows from \cite[Lemma 4.2]{S1} and \cite[Lemma 2.5]{NT} that if $h\leq k$, then $(\Omega_h,\phi_h, g_h, S_{\!\phi_h}(0,1))\in\mathcal{P}_{\lz,\,\Lz,\,\rz,\,C h^{1/2},*}$. It is easily seen that \cite[Lemma 2.7]{NT} (quadratic separation on $\partial U_h\cap B_c$) and \cite[Proposition 2.8]{NT} ($C^{2,\az}$ estimates in $U_h\cap B_c$) hold for the solution $w$ of \eqref{s3:NT2.71}. Hence we have the following version of \cite[Proposition 2.12]{NT} corresponding to our new definition of $\mathcal{P}_{\lz,\,\Lz,\,\rz,\,\kappa,\,\az}$.

\begin{prop}\label{s3:prop NT2.12}
Let $\Omega$ and $\phi$ satisfy \eqref{s3:1}-\eqref{s3:4}. Assume in addition that $\partial\Omega\cap B_\rz$ is $C^{2,\az}$ and $\phi\in C^{2,\az}(\partial\Omega\cap B_\rz)$ for some $\az\in (0,1)$. Then there exists $h_0>0$ depending only on $n,\lz,\Lz,\rz,\az,\|\partial\Omega\cap B_\rz\|_{C^{2,\az}}$ and $\|\phi\|_{C^{2,\az}(\partial\Omega\cap B_\rz)}$ such that for $h\leq h_0$, we have
\begin{eqnarray*}
\l(\Omega_h,\phi_h, g_h, S_{\!\phi_h}(0,1)\r)\in \mathcal{P}_{\lz,\,\Lz,\,\rz,\,C h^{1/2},\,\az}\quad
\mathrm{and}\quad\|\partial\Omega_h\cap B_{1/k}\|_{C^{2,\az}}\leq C'h^{1/2}.
\end{eqnarray*}
Here $k,C$ depend only on $n,\lz,\Lz$ and $\rz$. $C'$ depends only on $n,\lz,\Lz,\rz,\|\partial\Omega\cap B_\rz\|_{C^{2,\az}}$ and $\|\phi\|_{C^{2,\az}(\partial\Omega\cap B_\rz)}$.
\end{prop}

The covering lemma \cite[Lemma 3.13]{NT} holds for the newly defined class $\mathcal{P}_{\lz,\,\Lz,\,\rz,\,\kappa,\,\az}$.

\begin{lem}\label{s3:lem NT3.13}
Assume $(\Omega,\phi, g, U)\in \mathcal{P}_{\lz,\,\Lz,\,\rz,\,\kappa,*}$. Let $c$ be as in $(vi)$ in the definition of $\mathcal{P}_{\lz,\,\Lz,\,\rz,\,\kappa,*}$ and $w$ be the solution to \eqref{s3:P} for some constant $A>0$ such that $c(n,\lz,\Lz)\leq A\leq C(n,\lz,\Lz)$. Let $\psi$ denote one of the functions $\phi$ and $w$. Then there exists a sequence of disjoint sections $\{S_{\!\psi}(y_i,\dz_0\bar{h}(y_i))\}_{i=1}^\infty$, where $\dz_0=\dz_0(n,\lz,\Lz,\rz), y_i\in U\cap B_{c^2}$ and $S_{\!\psi}(y_i,\bar{h}(y_i))$ is the maximal interior section of $\psi$ in $U$, such that
\begin{equation}\label{s3:lem NT3.131}
U\cap B_{c^2}\subset\bigcup_{i=1}^\infty S_{\!\psi}\l(y_i,\f{\bar{h}(y_i)}{2}\r).
\end{equation}
Moreover, we have
\begin{equation}\label{s3:lem NT3.132}
S_{\!\psi}(y_i,\bar{h}(y_i))\subset U\cap B_c,\quad \bar{h}(y_i)\leq c.
\end{equation}
Let $M_d^{\mathrm{loc}}$ denote the number of sections $S_{\!\psi}(y_i,\bar{h}(y_i)/2)$ such that $d/2<\bar{h}(y_i)\leq d\leq c$, then
\begin{equation}\label{s3:lem NT3.133}
M_d^{\mathrm{loc}}\leq C_b d^{\f{1}{2}-\f{n}{2}}
\end{equation}
for some constant $C_b$ depending only on $n,\lz,\Lz,\rz$ and $\kappa$.
\end{lem}

Next we prove a stability result at the boundary for the cofactor matrices of $D^2\phi$. 

\begin{prop}\label{s3:prop NT3.14}
Given $0<\ez<1$. Assume $(\Omega,\phi, g, U)\in \mathcal{P}_{\lz,\,\Lz,\,\rz,\,\kappa,*}$. Let $c$ be as in $(vi)$ in the definition of $\mathcal{P}_{\lz,\,\Lz,\,\rz,\,\kappa,*}$ and $\tilde{c}$ be a constant depending only on $n,\lz,\Lz,\rz$ such that $\tilde{c}\leq c^2$. Assume $w$ is the solution to
$$\left\{
\begin{array}{rcl}
\mathrm{det}D^2 w=1&&{\mathrm{in}\;U},\\
w=\f{\phi}{(g^1)_{U\cap B_{\tilde{c}}}}&&{\mathrm{on}\;\partial U}.
\end{array}\right.$$
For any $1<p<\infty$, there exists $\ez=\ez(n,\lz,\Lz,\rz,p)>0$ such that if
$$\mathrm{mosc}_U g\leq\ez\quad\mathrm{and}\quad\sup_{S_{\!\phi}(x,h)\Subset U\cap B_c}\mathrm{mosc}_{S_{\!\phi}(x,h)}\,g\leq\ez,$$
then the following statements hold:
\begin{enumerate}
\item[\rm(i)]
\begin{equation*}\label{s3:prop NT3.141}
\|D^2\phi-(g^1)_{U\cap B_{\tilde{c}}}D^2 w\|_{L^p(U\cap B_{c^2})}\leq C(n,\lz,\Lz,\rz,\kappa, p)\ez^{\f{\dz}{n(2p-\dz)}}.
\end{equation*}
\item[\rm(ii)]\,Assume in addition that $(\Omega,\phi, g, U)\in \mathcal{P}_{\lz,\,\Lz,\,\rz,\,\kappa,\,\az}$.
Then
\begin{eqnarray*}
\l\|\Phi-(g^1)^{n-1}_{U\cap B_{\tilde{c}}}W \r\|_{L^p(U\cap B_{c^2})}
&\leq&C(n,\lz,\Lz,\rz,\kappa, p,\az)\ez^{\f{\dz}{n(2pn-\dz)}}.
\end{eqnarray*}
Here $\dz=\dz(n,\lz,\Lz,\rz)\in(0,1/2)$, $\Phi$ and $W$ are the cofactor matrices if $D^2\phi$ and $D^2 w$ respectively.
\end{enumerate}

\begin{proof}
The proof is similar to that of \cite[Proposition 3.14]{NT}, using Lemma \ref{s3:lem NT3.13} and Theorem \ref{s2:thm HA(i)} instead due to the $\mathrm{VMO}$-type condition of $g$. Note that if we take $A:=[g^1_{U\cap B_{\tilde{c}}}]^{-1}$, then $c(n,\lz,\Lz)\leq A\leq C(n,\lz,\Lz)$, as required in Lemma \ref{s3:lem NT3.13}. We sketch the proof.\\

\noindent(i) As in \cite{NT}, the conclusion follows from the two claims below and the interpolation inequality.\\

\noindent$\mathbf{Claim\,1.}$ For any $1\leq p<\infty$, there exists $C_0>0$ depending only on $n,\lz,\Lz,\rz,\kappa, p$ such that
$$\|D^2\phi\|_{L^p(U\cap B_{c^2})}+\|D^2 w\|_{L^p(U\cap B_{c^2})}\leq C_0$$
if $\ez=\ez(n,\lz,\Lz,\rz,p)>0$ is small.

\noindent$\mathbf{Claim\,2.}$ There exist $\dz=\dz(n,\lz,\Lz,\rz)\in(0,1/2)$ and $C=C(n,\lz,\Lz,\rz,\kappa)$ such that
\begin{equation*}
\|D^2\phi-(g^1)_{U\cap B_{\tilde{c}}}D^2 w\|_{L^\dz(U\cap B_{c^2})}\leq C\ez^{1/n}
\end{equation*}
for any $0<\ez<1$.

Let $\psi$ denote one of the functions $\phi$ and $w$. By Lemma \ref{s3:lem NT3.13}, $\mathbf{Claim\,1}$ follows from this result: there exists $C=C(n,\lz,\Lz,\rz,p)>0$ such that
\begin{equation}\label{s3:prop NT3.141}
\int_{S_{\!\psi}\l(y,\f{\bar{h}(y)}{2}\r)}|D^2\psi|^p\leq C\bar{h}(y)^{\f{n}{2}}|\log\bar{h}(y)|^{2p}\quad\forall y\in U\cap B_{c^2}.
\end{equation}
To prove the above inequality, we denote $h:=\bar{h}(y)\leq c$. Apply Proposition \ref{s3:prop NT2.3} to $S_{\!\psi}(y,h)$ and we obtain
$$k_0 E_h\subset S_{\!\psi}(y,h)-y\subset k_0^{-1}E_h,$$
where $E_h:=h^{\f{1}{2}}A_h^{-1}B_1$ with $\mathrm{det}\,A_h=1$ and $\|A_h\|,\|A_h^{-1}\|\leq C|\log h|$. Define 
$$\tilde{\psi}_h(x):=\f{1}{h}[\psi(y+h^{\f{1}{2}}A_h^{-1}x)-l_{\psi,y}(h^{\f{1}{2}}A_h^{-1}x)-h],\quad x\in\tilde{\Omega}_h:=h^{-\f{1}{2}}A_h(\Omega-y)$$
If $\psi=\phi$, we denote $\tilde{g}_h(x):=g(y+h^{\f{1}{2}}A_h^{-1}x)$. Denote $T_hx:=h^{-\f{1}{2}}A_h(x-y)$, then
\begin{eqnarray*}
\sup_{S_{\!\tilde{\phi}_h}(x,t)\Subset S_{\!\tilde{\phi}_h}(0,1)}\mathrm{mosc}_{S_{\!\tilde{\phi}_h}(x,t)}\,\tilde{g}_h
&\leq&\sup_{S_{\!\phi}(T_h^{-1}x,\,t h)\Subset S_{\!\phi}(y,h)}\mathrm{mosc}_{S_{\!\phi}(T_h^{-1}x,\,th)}\,g\\
&\leq&\sup_{S_{\!\phi}(T_h^{-1}x,\,th)\Subset U\cap B_c}\mathrm{mosc}_{S_{\!\phi}(T_h^{-1}x,\,th)}\,g
\leq\ez.
\end{eqnarray*}
Therefore by Theorem \ref{s2:thm HA(i)} $(i)$ we have
$$\int_{S_{\!\tilde{\phi}_h}\l(0,\f{1}{2}\r)}\|D^2\tilde{\phi}_h\|^p\leq C(n,\lz,\Lz,\rz,p)$$
if $\ez=\ez(n,\lz,\Lz,\rz,p)$ is small.

If $\psi=w$ then a similar inequality holds for $w$ because of Pogorelov's estimate \cite[(4.2.6)]{G}. Changing variables as in Page 662 in \cite{NT}, we obtain \eqref{s3:prop NT3.141} and then $\mathbf{Claim\,1}$ follows.

The proof of $\mathbf{Claim\,2}$ is similar to that of Lemma \ref{s2:lem GN13.4}. Let $v:=\phi-(g^1)_{U\cap B_{\tilde{c}}} w$, then $v=0$ on $\partial U$. Note that $|U\cap B_{\tilde{c}}|\geq c(n,\lz,\Lz,\rz)$, then by Lemma \ref{s2:lem H3.1} and Proposition \ref{s2:prop VMO},
\begin{eqnarray*}\label{s3:prop NT3.142}
\max_{\overline{U}}|\phi-(g^1)_{U\cap B_{\tilde{c}}} w|
&\leq&C_n\mathrm{diam}(U)\left(\int_{U}|g^1-(g^1)_{U\cap B_{\tilde{c}}}|^{n}dx\right)^{\f{1}{n}}\nonumber\\
&\leq&C(n,\lz,\Lz,\rz)\left(\int_{U}|g^1-(g^1)_{U}|^{n}dx\right)^{\f{1}{n}}\nonumber\\
&\leq&C(n,\lz,\Lz,\rz)\ez^{\f{1}{n}}.
\end{eqnarray*}
Using this and similar arguments to the proof of \cite[Proposition 3.14]{NT}, $\mathbf{Claim\;2}$ is proved and the conclusion (i) follows.\\

\noindent(ii) Write
$$\Phi-(g^1)^{n-1}_{U\cap B_{\tilde{c}}}W=\l[1-\f{(g^1)^{n}_{U\cap B_{\tilde{c}}}}{\mathrm{det}D^2\phi}\r]\Phi
-\f{(g^1)^{n-1}_{U\cap B_{\tilde{c}}}}{\mathrm{det}D^2\phi}\Phi\l(D^2\phi-(g^1)^{n-1}_{U\cap B_{\tilde{c}}}D^2 w\r)W.$$
Arguing as in the proof of Lemma \ref{s2:lem GN13.4}, for any $1\leq q,r<\infty$ we have,
\begin{eqnarray*}
&&\l\|\Phi-(g^1)^{n-1}_{U\cap B_{\tilde{c}}}W \r\|_{L^q(U\cap B_{c^2})}\\
&\leq&C(n,\lz,\Lz,\rz,q,r)\l(\ez^{\f{1}{n}}+\ez^{\f{1}{qr}}
+\|D^2 w\|^{n-1}_{L^\infty(U\cap B_{c^2})}\l\|D^2\phi-(g^1)^{n-1}_{U\cap B_{\tilde{c}}}D^2 w\r\|_{L^{qr}(U\cap B_{c^2})}\r)\|D^2\phi\|^{n-1}_{L^{qr'(n-1)}(U\cap B_{c^2})}.
\end{eqnarray*}
Choose $r=n$ and then $r'=\f{n}{n-1}$, if $\ez=\ez(n,\lz,\Lz,\rz,q)$ is small, then
\begin{eqnarray*}
\l\|\Phi-(g^1)^{n-1}_{U\cap B_{\tilde{c}}}W \r\|_{L^q(U\cap B_{c^2})}\leq C(n,\lz,\Lz,\rz,q,\kappa)\l(\ez^{\f{1}{n}}+\ez^{\f{1}{qn}}
+\|D^2 w\|^{n-1}_{L^\infty(U\cap B_{c^2})}\ez^{\f{\dz}{n(2qn-\dz)}}\r)
\end{eqnarray*}
by (i) and $\mathbf{Claim\,1}$ in the proof of (i). By the definition of $\mathcal{P}_{\lz,\,\Lz,\,\rz,\,\kappa,\az}$, we have
$$\|D^2 w\|^{n-1}_{L^\infty(U\cap B_{c^2})}\leq C(n,\lz,\Lz,\rz,\az).$$
Hence the conclusion follows.
\end{proof}
\end{prop}

We next list the global version of Lemma \ref{s2:lem GN22.7}. Recall Section $2$ for the definition of $d(x,x_0)$ and $G_{\!M}(u,\Omega)$. Assume 
\begin{equation}\label{s3:g1}
\Omega\subset B_{1/\rz}\;\mathrm{and\;for\;each}\;y\in\partial\Omega\;\mathrm{there\;is\;a\;ball\;}B_\rz(z)\subset\Omega\;
\mathrm{that\;is\;tangent\;to\;}\partial\Omega\;\mathrm{at}\;y.
\end{equation}
Let $\phi:\overline{\Omega}\rightarrow\mathbb{R}$, $\phi\in C^{0,1}(\overline{\Omega})\cap C^2(\Omega)$ be a convex function satisfying
\begin{eqnarray}\label{s3:g3}
\mathrm{det}D^2\phi=g,\quad\quad 0<\lambda\leq g\leq\Lambda\quad\mathrm{in}\;\Omega.
\end{eqnarray}
Assume further that on $\partial\Omega$, $\phi$ separates quadratically from its tangent planes, namely,
\begin{equation}\label{s3:g4}
\rz|x-x_0|^2\leq\phi(x)-\phi(x_0)-\nabla\phi(x_0)\cdot (x-x_0)\leq\rz^{-1}|x-x_0|^2,\quad\quad\forall x,x_0\in\partial\Omega.
\end{equation}

\begin{lem}(See \cite[Lemma 3.4]{NT}.)\label{s3:lem NT3.4}
Assume $\Omega$ satisfies \eqref{s3:g1} and $\phi\in C(\overline{\Omega})$ is a solution of \eqref{s1:2}. Define
$$A_\sz^{\mathrm{loc}}:=\{x_0\in\Omega: d(x,x_0)^2\geq\sz|x-x_0|^2,\;\mathrm{for\;all}\;x\;\mathrm{in\;some\;neighborhood\;of\;}x_0\}.$$
Let $u\in C^2(\Omega)$. Then for $\kappa>1$, we have
\begin{equation}\label{s3:NT3.41}
\{x\in\Omega:|D_{i j}u(x)|>\bz^\kappa\}\subset\l(\Omega\backslash A^{\mathrm{loc}}_{(c\bz^{\f{\kappa-1}{2}})^{\f{-2}{n-1}}}\r)\cup\l(\Omega\backslash G_{\bz}(u,\Omega)\r)
\end{equation}
for any $\bz>0$, where $c>0$ is a universal constant.
\end{lem}

The estimate of the first term on the right-hand side of \eqref{s3:NT3.41} is given by the following theorem.

\begin{thm}\label{s3:thm NT3.5}
Assume $\Omega, \phi$ satisfy \eqref{s3:g1}-\eqref{s3:g4} and $\partial\Omega\in C^{1,1}$. Let $0<\ez<1$. Suppose $$\sup_{S_{\!\phi}(x,h)\Subset\Omega}\mathrm{mosc}_{S_{\!\phi}(x,h)}g\leq\ez.$$
Then there exist universal constants $M>0,0<\dz_1<1$ such that
\begin{equation*}\label{s3:NT3.51}
|\Omega\backslash A^{\mathrm{loc}}_{s^{-2}}|\leq C'(\ez,n,\lz,\Lz,\rz,\|\partial\Omega\|_{C^{1,1}})s^{\f{\ln\l(2\ez^{\dz_1/2}\r)}{\ln M}}\quad\forall s>0.
\end{equation*}
In particular, for $s=(c\bz^{\f{\kappa-1}{2}})^{\f{1}{n-1}}$ we have
\begin{equation*}\label{s3:NT3.52}
|\Omega\backslash A^{\mathrm{loc}}_{(c\bz^{\f{\kappa-1}{2}})^{\f{-2}{n-1}}}|\leq C'(\ez,n,\lz,\Lz,\rz,\|\partial\Omega\|_{C^{1,1}})\,\bz^{\f{(\kappa-1)\ln\l(2\ez^{\dz_1/2}\r)}{2(n-1)\ln M}}\quad\forall\,\bz>0.
\end{equation*}

\begin{proof}
The proof is similar to that of \cite[Theorem 3.5]{NT}, using Lemma \ref{s2:lem H6.3} instead. Namely, we cover $\Omega$ by sections of $\phi$ given by \cite[Lemma 3.12]{NT} and obtain 
$$|\Omega\backslash A^{\mathrm{loc}}_{s^{-2}}|\le\sum_{h_i\leq c}\l|S_{\!\phi}\l(y_i,h_i/2\r)\backslash A^{\mathrm{loc}}_{s^{-2}}\r|+\sum_{h_i> c}\l|S_{\!\phi}\l(y_i,h_i/2\r)\backslash A^{\mathrm{loc}}_{s^{-2}}\r|:=I+II,$$
where $h_i:=\bar{h}(y_i)$. For the summation $I$ corresponding to $h\leq c$, let $Tx:=h^{-\f{1}{2}}A_{h}(x-y)$ where $A_{h}$ is given by Proposition \ref{s3:prop NT2.3}. Define the rescaled domain $\tilde{U}_h:=T(S_{\!\phi}(y,h))$ and function 
$$\tilde{\phi}_h(x):=\f{1}{h}[\phi(T^{-1}x)-l_{\phi,y}(T^{-1}y)-h].$$
Denote $\tilde{g}_h=\mathrm{det}\,D^2\tilde{\phi}_h$ and
$$\tilde{D}^{\f{1}{2}}_s:=\{x\in S_{\!\tilde{\phi}_h}(0,1/2): S_{\!\tilde{\phi}_h}(x,t)\subset B(x,s\sqrt{t}),\quad\forall t\le\eta\}.$$
We apply Lemma \ref{s2:lem H6.3} with $\Omega\rightsquigarrow\tilde{U}_h, \phi\rightsquigarrow\tilde{\phi}_h, g\rightsquigarrow\tilde{g}_h, \az\rightsquigarrow \f{1}{2}$. Denote $\az_0:=\f{\az+1}{2}=\f{3}{4}$. Note that for any $S_{\!\tilde{\phi}_h}(x,t)\subset(\tilde{U}_h)_{\f{\az_0+1}{2}}\Subset\tilde{U}_h$, we have
$$\mathrm{mosc}_{S_{\!\tilde{\phi}_h}(x,t)}\tilde{g}_h=\mathrm{mosc}_{S_{\!\phi}(T^{-1}x,\,th)}g\leq\ez.$$
Thus by Lemma \ref{s2:lem H6.3}, there exist universal constants $M>0,0<\dz_1<1$ such that
\begin{equation*}\label{s3:NT3.53}
|S_{\!\tilde{\phi}_h}(0,1/2)\backslash \tilde{D}^{\f{1}{2}}_s|\leq\f{|\tilde{U}_h|}{(2\ez^{\dz_1/2})^2}s^{\f{\ln\l(2\ez^{\dz_1/2}\r)}{\ln M}}.
\end{equation*}
Then we obtain 
$$I\le C(n,\lz,\Lz,\rz,\ez,\|\partial\Omega\|_{C^{1,1}})s^{-p_\ez},\quad p_{\ez}=-\f{\ln\l(2\ez^{\dz_1/2}\r)}{\ln M}.$$
The estimate of summation $II$ is also similar to that in \cite{NT} by using standard normalization for interior sections. Again, we use Lemma \ref{s2:lem H6.3} in each of these sections.
\end{proof}
\end{thm}

\section{Boundary estimate-Power decay estimates for the solutions of \eqref{s1:1}}\label{s6}

In this section, we prove a localized estimate at the boundary of the density of good sets where $D^2\phi$ is bounded. Different from \cite[Lemma 4.6]{NT}, since $g$ is in $\mathrm{VMO}$-type spaces, there is a difficulty in using convex envelopes of functions to help estimate the density of good sets. Instead, we compare $\phi$ with the solution $w$ of $\mathrm{det}\,D^2w=1$. We apply the localized covering lemma at the boundary and cover $U\cap B_{c^2}$ by sections of $w$. Then we use the one-sided small power decay estimate \cite[Proposition 3.4]{GN2} in each of these sections to the linearized operator $\mathcal{L}_w$ to estimate the set where $\phi$ can be touched from above by a quadratic polynomial. The precise statement is as follows.

\begin{lem}\label{s3:lem NT4.6}
Assume $(\Omega,\phi, g, U)\in \mathcal{P}_{\lz,\,\Lz,\,\rz,\,\kappa,\az}$ with $0<\kappa<1$, $\phi\in C^1(\Omega)$ and $\mathrm{mosc}_U g\leq\ez$. Let $c$ be the constant in the definition of the class $\mathcal{P}$. Then there exist $\sz=\sz(n,\lz,\Lz,\rz,\az,\kappa)>0$, $0<\tilde{c}=\tilde{c}(n,\lz,\Lz,\rz)\leq c^4$ and $0<\tau'=\tau'(n)<1/2$ such that
\begin{eqnarray*}
|(U\cap B_{\tilde{c}})\backslash A_{\sz}|\leq\ez^{\f{\tau'}{n}}|U\cap B_{\tilde{c}}|,
\end{eqnarray*}
and for any $S_{\!\phi}(0,r)$ with $r=r(n,\lz,\Lz,\rz)\leq\tilde{c}^3$, we have
\begin{eqnarray*}
|S_{\!\phi}(0,r)\backslash A_{\sz}|\leq\ez^{\f{\tau'}{n}}|S_{\!\phi}(0,r)|,
\end{eqnarray*}
where
\begin{equation}\label{s3:NT4.61}
A_\sz:=\{\bar{x}\in U\cap B_{\tilde{c}}:\phi(x)\geq l_{\phi,\bar{x}}(x)+\sz|x-\bar{x}|^2,\quad\forall x\in U\cap B_{\tilde{c}}\}.
\end{equation}

\begin{proof}
Let $w$ be the convex solution of 
$$\left\{
\begin{array}{rcl}
\mathrm{det}D^2 w=1&&{\mathrm{in}\;U},\\
w=\f{\phi}{(g^1)_{U}}&&{\mathrm{on}\;\partial U}.
\end{array}\right.$$
By Lemma \ref{s2:lem H3.1} and Proposition \ref{s2:prop VMO},
\begin{eqnarray}\label{s3:NT4.62}
\max_{\overline{U}}|\phi-(g^1)_{U} w|
&\leq&C_n\mathrm{diam}(U)\left(\int_{U}|g^1-(g^1)_{U}|^{n}dx\right)^{\f{1}{n}}\nonumber\\
&\leq&C(n,\lz,\Lz,\rz)\ez^{\f{1}{n}}.
\end{eqnarray}
Consider the operator $\mathcal{M}u:=(\mathrm{det}\,D^2 u)^{\f{1}{n}}$ and its linearized operator
$$\hat{L}_u v:=\frac{1}{n}(\mathrm{det}D^2 u)^{\f{1}{n}}\mathrm{tr}[(D^2 u)^{-1}D^2 v],$$
Denote $v:=\f{\phi}{(g^1)_U}-w$. Since $\mathcal{M}$ is concave, we have
\begin{eqnarray}\label{s3:NT4.63}
F:=\f{g^1}{(g^1)_{U}}-1=\mathcal{M}\l(\f{\phi}{(g^1)_U}\r)-\mathcal{M} w\leq \hat{L}_{w} v=\f{1}{n}\mathrm{tr}[(D^2 w)^{-1}D^2 v]\quad\quad\mathrm{in}\;U.
\end{eqnarray}
Let $c$ be the constant in the definition of the class $\mathcal{P}$. For $u\in C^1(U\cap B_c)$ and $\bz>0$ we define
\begin{eqnarray}\label{s3:NT4.64}
G^+_\bz(u,w):=\l\{\bar{x}\in U\cap B_c: u(x)-l_{u,\bar{x}}(x)\leq\f{\bz}{2}d_w(x,\bar{x})^2,\quad\forall x\in U\cap B_c\r\},
\end{eqnarray}
where $d_w(x,\bar{x})^2:=w(x)-l_{w,\bar{x}}(x)$ for any $x,\bar{x}\in U\cap B_c$.

By Lemma \ref{s3:lem NT3.13}, there exists a sequence of disjoint sections $\{S_{\!w}(y_i,\dz_0\bar{h}(y_i))\}_{i=1}^\infty$, where $\dz_0=\dz_0(n,\lz,\Lz,\rz), y_i\in U\cap B_{c^2}$ and $S_w(y_i,\bar{h}(y_i))$ is the maximal interior section of $w$ in $U$, such that
\begin{equation*}
U\cap B_{c^2}\subset\bigcup_{i=1}^\infty S_{\!w}\l(y_i,\f{\bar{h}(y_i)}{2}\r)
\end{equation*}
and
\begin{equation*}
S_{\!w}(y_i,\bar{h}(y_i))\subset U\cap B_c,\quad \bar{h}(y_i)\leq c.
\end{equation*}
Moreover, let $M_d^{\mathrm{loc}}$ denote the number of sections $S_{\!w}(y_i,\bar{h}(y_i)/2)$ such that $d/2<\bar{h}(y_i)\leq d\leq c$, then
\begin{equation}\label{s3:NT4.65}
M_d^{\mathrm{loc}}\leq C_b d^{\f{1}{2}-\f{n}{2}}
\end{equation}
for some constant $C_b$ depending only on $n,\lz,\Lz,\rz$ and $\kappa$.

Let $h:=\bar{h}(y)$ with $\bar{h}(y)\leq c$. Let $T$ be an affine map such that
\begin{equation*}\label{s3:NT4.66}
B_{\az_n}\subset T(S_{\!w}(y,h))\subset B_1.
\end{equation*}

Let $\tilde{U}:=T(S_{\!w}(y,h))$ and $\tilde{\Omega}:=T(U\cap B_c)$. For each $\tilde{x}\in\tilde{\Omega}$, define
$$\tilde{w}(\tilde{x})=|\mathrm{det}\,T|^{\f{2}{n}}\l[w(T^{-1}\tilde{x})-l_{w,y}(T^{-1}\tilde{x}) -h\r],\quad\mathrm{and}\quad\tilde{v}(\tilde{x})=v(T^{-1}\tilde{x}).$$
Then by \eqref{s3:NT4.63} we have
$$\mathrm{det}D^2\tilde{w}=1\;\mathrm{in}\;\tilde{\Omega},
\quad\quad\tilde{w}=0\;\mathrm{on}\;\partial\tilde{U},$$
\begin{eqnarray*}
\mathcal{L}_{\tilde{w}}\tilde{v}(\tilde{x})\geq\tilde{F}(\tilde{x}):=n|\mathrm{det}\,T|^{\f{-2}{n}}F(T^{-1}\tilde{x})\quad\mathrm{in}\;\tilde{\Omega}.
\end{eqnarray*}

Applying \cite[Proposition 3.4]{GN2} (a straightforward modification of the proof gives rise to one-sided estimates) with $U\rightsquigarrow\tilde{U},\Omega\rightsquigarrow\tilde{\Omega},\phi\rightsquigarrow\tilde{w}, u\rightsquigarrow\tilde{v},\az\rightsquigarrow \f{1}{2}$, we obtain
\begin{equation}\label{s3:NT4.67}
\l|\tilde{U}_{\f{1}{2}}\backslash G^+_\bz(\tilde{v},\tilde{\Omega},\tilde{w})\r|
\leq C\f{[\|\tilde{v}\|_{L^\infty(\tilde{\Omega})}+\|\tilde{F}\|_{L^n(\tilde{U})}]^{\tau'}}{\bz^{\tau'}}\quad\forall\bz>0,
\end{equation}
where $C,\tau'>0$ depending only on $n$ and
$$G^+_\bz(\tilde{v},\tilde{\Omega},\tilde{w}):=\l\{\bar{x}\in\tilde{\Omega}:\tilde{v}(\tilde{x})-l_{\tilde{v},\bar{x}}(\tilde{x})
\leq\f{\bz}{2}[\tilde{w}(\tilde{x})-l_{\tilde{w},\bar{x}}(\bar{x})],\quad\forall \tilde{x}\in\tilde{\Omega}\r\}.$$
Note that $\tilde{U}_{\f{1}{2}}=T\l(S_{w}\l(y,\f{h}{2}\r)\r)$ and
\begin{eqnarray*}
G^+_\bz(\tilde{v},\tilde{\Omega},\tilde{w})=T(G^+_{\bz|\mathrm{det}\,T|^{\f{2}{n}}}(v,w)),
\end{eqnarray*}
where we recall \eqref{s3:NT4.64} for the definition of $G^+_{\bz|\mathrm{det}\,T|^{\f{2}{n}}}(v,w)$. Moreover,
$$\|\tilde{F}\|_{L^n(\tilde{U})}=\f{n}{|\mathrm{det}\,T|^{\f{1}{n}}}\|F\|_{L^n(S_w(y,h))}.$$
Thus \eqref{s3:NT4.67} implies that
\begin{equation}\label{s3:NT4.68}
\l|S_{\!w}\l(y,\f{h}{2}\r)\backslash G^+_{\bz|\mathrm{det}\,T|^{\f{2}{n}}}(v,w)\r|
\leq \f{C}{\bz^{\tau'}}|\mathrm{det}\,T|^{-1}\l[\|v\|_{L^\infty(U\cap B_c)}+\f{n}{|\mathrm{det}\,T|^{\f{1}{n}}}\|F\|_{L^n(S_{\!w}(y,h))}\r]^{\tau'}\quad\forall\bz>0.
\end{equation}
By $(vii)$ in the definition of $\mathcal{P}$,
\begin{eqnarray}\label{s3:NT4.69}
\f{c_0}{2}|x-x_0|^2\le w(x)-l_{w,x_0}(x)\leq\f{c_0^{-1}}{2}|x-x_0|^2\quad\forall x,x_0\in U\cap B_c,
\end{eqnarray}
where $c_0=c_0(n,\lz,\Lz,\rz,\az)$. Hence for any $\bar{x}\in G^+_{\bz|\mathrm{det}\,T|^{\f{2}{n}}}(v,w)$ and $x\in U\cap B_c$, we have 
$$v(x)-l_{v,\bar{x}}(x)\leq\f{C_0\,\bz h^{-1}}{2}|x-\bar{x}|^2,$$
where $C_0=C_0(n,\lz,\Lz,\rz,\az)$. Therefore, by \eqref{s3:NT4.68} and \eqref{s3:NT4.62},
\begin{eqnarray*}\label{s3:NT4.611}
&&\l|S_{\!w}\l(y,\f{h}{2}\r)\backslash\{\bar{x}\in U\cap B_c: v(x)-l_{v,\bar{x}}(x)
\leq\f{\bz'}{2}|x-\bar{x}|^2,\quad\forall x\in U\cap B_c\}\r|\nonumber\\
&\leq&\f{C(n,\lz,\Lz,\rz,\az)}{\bz'^{\tau'}}h^{\f{n}{2}-\tau'}\ez^{\f{\tau'}{n}},\quad\forall\bz'>0,
\end{eqnarray*}
which implies
\begin{eqnarray*}\label{s3:NT4.612}
&&\l|(U\cap B_{c^2})\backslash \{\bar{x}\in U\cap B_c: v(x)-l_{v,\bar{x}}(x)
\leq\f{\bz'}{2}|x-\bar{x}|^2,\quad\forall x\in U\cap B_c\}\r|\nonumber\\
&\leq&\f{C\ez^{\f{\tau'}{n}}}{\bz'^{\tau'}}\sum_{k=0}^\infty\sum_{i\in\mathcal{F}_{c 2^{-k}}}\bar{h}(y_i)^{\f{n}{2}-\tau'}\nonumber\\
&\leq&\f{C\ez^{\f{\tau'}{n}}}{\bz'^{\tau'}}\sum_{k=0}^\infty(c 2^{-k})^{\f{n}{2}-\tau'}M^{\mathrm{loc}}_{c 2^{-k}}\leq\f{C\ez^{\f{\tau'}{n}}}{\bz'^{\tau'}},\quad\forall\bz'>0.
\end{eqnarray*}
if we choose $0<\tau'<\f{1}{2}$, where we use \eqref{s3:NT4.65} and $C=C(n,\lz,\Lz,\rz,\az,\kappa)$.

For $x,\bar{x}\in U\cap B_c$ such that $v(x)-l_{v,\bar{x}}(x)\leq\f{\bz'}{2}|x-\bar{x}|^2$, we have by the definition of $v$ and \eqref{s3:NT4.69}
\begin{eqnarray*}
\f{1}{(g^1)_U}[\phi(x)-l_{\phi,\bar{x}}(x)]&\leq&\f{\bz'}{2}|x-\bar{x}|^2
+w(x)-l_{w,\bar{x}}(x)\leq\f{\bz'+c_0^{-1}}{2}|x-\bar{x}|^2.
\end{eqnarray*}
It follows that
\begin{eqnarray}\label{s3:NT4.613}
&&\l|(U\cap B_{c^2})\backslash \{\bar{x}\in U\cap B_c: \phi(x)-l_{\phi,\bar{x}}(x)\leq\Lz^{\f{1}{n}}\f{\bz'+c_0^{-1}}{2}|x-\bar{x}|^2,\quad\forall x\in U\cap B_c\}\r|\nonumber\\
&\leq&\f{C\ez^{\f{\tau'}{n}}}{\bz'^{\tau'}}\leq\f{C(n,\lz,\Lz,\rz,\az,\kappa)\ez^{\f{\tau'}{n}}}{\bz'^{\tau'}}|U\cap B_{c^2}|,\quad\forall\bz'>0.
\end{eqnarray}

For any fixed $\bz'>0$, we denote $M_1=M_1(\bz'):=\Lz^{\f{1}{n}}\f{\bz'+c_0^{-1}}{2}$. Let $0<\tilde{c}\leq c^2$ be a constant to be chosen later. Fix $\bar{x}\in U\cap B_{\tilde{c}}$ such that $\phi(x)-l_{\phi,\bar{x}}(x)\leq M_1|x-\bar{x}|^2,\;\forall x\in U\cap B_{c}$. Then we have
\begin{eqnarray}\label{s3:NT4.614}
(U\cap B_c)\cap B_{\sqrt{\f{h}{M_1}}}(\bar{x})\subset S_{\!\phi}(\bar{x},h),\quad\forall h>0.
\end{eqnarray}

Fix any $h\leq c^4$, then we have $B_{\sqrt{\f{h}{M_1}}}(\bar{x})\subset B_{2 c^2}\subset B_{c/2}$. Let $Tx=Ax+b$ be an affine transformation such that
\begin{equation}\label{s3:NT4.623}
B_{\az_n}\subset T(S_{\!\phi}(\bar{x},h))\subset B_1.
\end{equation}
By the volume estimate for boundary sections we have
\begin{equation}\label{s3:NT4.624}
|\mathrm{det}\,T|\approx C(n,\lz,\Lz,\rz)h^{-\f{n}{2}}.
\end{equation}

Let $x_0\in\partial U\cap B_{c}$ such that $|\bar{x}-x_0|=\mathrm{dist}(\bar{x},\partial U)$. Denote $\bar{y}:=\bar{x}+\f{1}{2}\sqrt{\f{h}{M_1}}\nu_{x_0}$, then $B_{\f{1}{2}\sqrt{\f{h}{M_1}}}(\bar{y})\subset B_{\sqrt{\f{h}{M_1}}}(\bar{x})$. We now prove that
\begin{equation}\label{s3:NT4.619}
\partial B_{\f{1}{2}\sqrt{\f{h}{M_1}}}(\bar{y})\cap\big((\partial U\cap B_c)\backslash\{x_0\}\big)=\emptyset.
\end{equation}
In fact, for any $x\in\partial B_{\f{1}{2}\sqrt{\f{h}{M_1}}}(\bar{y})$ we have 
\begin{equation}\label{s3:NT4.620}
(x-x_0)\cdot\nu_{x_0}\geq\sqrt{\f{h}{M_1}}[(x-x_0)\cdot\nu_{x_0}-|\bar{x}-x_0|]\geq\l|(x-x_0)-[(x-x_0)\cdot\nu_{x_0}]\nu_{x_0}\r|^2.
\end{equation}
This together with $(vi)$ in the definition of $\mathcal{P}$ implies \eqref{s3:NT4.619} and it follows that
\begin{equation*}
B_{\f{1}{2}\sqrt{\f{h}{M_1}}}(\bar{y})\subset(U\cap B_c).
\end{equation*}
This together with \eqref{s3:NT4.614} yields
\begin{equation}\label{s3:NT4.622}
B_{\f{1}{2}\sqrt{\f{h}{M_1}}}(\bar{y})\subset S_{\!\phi}(\bar{x},h).
\end{equation}

By \eqref{s3:NT4.622} and the second inclusion in \eqref{s3:NT4.623} we have $\|A\|\leq 4\sqrt{\f{M_1}{h}}$. This together with \eqref{s3:NT4.624} implies $\|A^{-1}\|\leq C\sqrt{M_1^{n-1}h}$, where $C=C(n,\lz,\Lz,\rz)$. Using the second conclusion in \eqref{s3:NT4.623} we obtain $\mathrm{diam}(S_{\!\phi}(\bar{x},h))\leq C\sqrt{M_1^{n-1}h}$. Therefore,
\begin{equation}\label{s3:NT4.625}
S_{\!\phi}(\bar{x},h)\subset B_{C\sqrt{M_1^{n-1}h}}(\bar{x}).
\end{equation}

For any $x\in U\cap B_{\tilde{c}}$, by $(v)$ in the definition of $\mathcal{P}$, there is a constant $C_0=C_0(n,\lz,\Lz,\rz)>0$ such that 
$$h:=\phi(x)-l_{\phi,\bar{x}}(x)\leq 2 C_0|x-\bar{x}|\leq 4 C_0\tilde{c}\leq c^4,$$
where we choose $\tilde{c}:=\f{c^4}{4 C_0}$. Thus \eqref{s3:NT4.625} gives
\begin{equation}\label{s3:NT4.626}
\phi(x)-l_{\phi,\bar{x}}\geq\f{1}{C^2 M_1^{n-1}}|x-\bar{x}|^2.
\end{equation}

Since $0<\tilde{c}\leq c^2$, the estimate \eqref{s3:NT4.613} also holds if we replace $U\cap B_{c^2}$ by $U\cap B_{\tilde{c}}$, hence by \eqref{s3:NT4.626}
\begin{eqnarray*}\label{s3:NT4.627}
&&\l|(U\cap B_{\tilde{c}})\backslash \{\bar{x}\in U\cap B_{\tilde{c}}: \phi(x)-l_{\phi,\bar{x}}\geq\f{1}{C^2 M_1^{n-1}}|x-\bar{x}|^2,\quad\forall x\in U\cap B_{\tilde{c}}\}\r|\nonumber\\
&\leq&\f{C(n,\lz,\Lz,\rz,\az,\kappa)\ez^{\f{\tau'}{n}}}{\bz'^{\tau'}}|U\cap B_{\tilde{c}}|,\quad\forall\bz'>0,
\end{eqnarray*}
and for any $S_{\!\phi}(0,r)$ with $r=r(n,\lz,\Lz,\rz)\leq\tilde{c}^3$,
\begin{eqnarray*}\label{s3:NT4.628}
&&\l|S_{\!\phi}(0,r)\backslash \{\bar{x}\in U\cap B_{\tilde{c}}: \phi(x)-l_{\phi,\bar{x}}(x)\geq\f{1}{C^2 M_1^{n-1}}|x-\bar{x}|^2,\quad\forall x\in U\cap B_{\tilde{c}}\}\r|\nonumber\\
&\leq&\f{C(n,\lz,\Lz,\rz,\az,\kappa)\ez^{\f{\tau'}{n}}}{\bz'^{\tau'}}|S_{\!\phi}(0,r)|,\quad\forall\bz'>0,
\end{eqnarray*}
where $M_1=M_1(\bz'):=\Lz^{\f{1}{n}}\f{\bz'+c_0^{-1}}{2}$ and $C=C(n,\lz,\Lz,\rz)$. Fix $\bz'$ large such that $\f{C(n,\lz,\Lz,\rz,\az,\kappa)}{\bz'^{\tau'}}\leq 1$ and denote $\sz:=\f{1}{C^2 M_1^{n-1}}$, and we obtain the conclusion.
\end{proof}
\end{lem}

The rest part of this section is devoted to the estimate of the second term on the right-hand side of \eqref{s3:NT3.41}. Firstly, a straightforward modification of \cite[Lemma 4.5]{NT} gives the lemma below. 

\begin{lem}\label{s3:lem NT4.5}
Assume $(\Omega,\phi, g, U)\in \mathcal{P}_{\lz,\,\Lz,\,\rz,\,\kappa,\az}$. Let $c$ be the constant in the definition of the class $\mathcal{P}$ and $\tilde{c}\leq c^2/2$ be the constant given by Lemma \ref{s3:lem NT4.6}. Suppose $u\in C(\overline{U})\cap W^{2,n}_{\mathrm{loc}}(U)$ is a solution of $\mathcal{L}_\phi u=f$ in $U\cap B_{c^2}$ with
$$\|u\|_{L^\infty(U\cap B_{c^2})}+\|u\|_{C^{2,\az}(\partial U\cap B_{c^2})}\leq 1.$$
Let $w$ be defined as in $(vii)$ of the definition of $\mathcal{P}$ with $A:=[(g^1)_{U\cap B_{\tilde{c}}}]^{-1}$ and $W$ be the cofactor matrix of $D^2 w$. Assume $h$ is a solution of
$$\left\{
\begin{array}{rcl}
\mathcal{L}_wh=0&&{\mathrm{in}\;U\cap B_{\tilde{c}}},\\
h=u&&{\mathrm{on}\;\partial(U\cap B_{\tilde{c}})}.
\end{array}\right.$$
Then there exist $C>0$ and $0<\gz<1$ depending only on $n,\lz,\Lz,\rz,\az$ such that
\begin{equation}\label{s3:NT4.51}
\|h\|_{C^{1,1}\l(\overline{U\cap B_{\tilde{c}/2}}\r)}\leq C,
\end{equation}
and if $\|\Phi-(g^1)^{n-1}_{U\cap B_{\tilde{c}}W}\|_{L^n(U\cap B_{\tilde{c}})}\leq\l(\f{\tilde{c}}{2}\r)^{4}$ then
\begin{eqnarray*}
\|u-h\|_{L^\infty\l(U\cap B_{\tilde{c}/2}\r)}&+&\l\|f-\mathrm{tr}\l([\Phi-(g^1)^{n-1}_{U\cap B_{\tilde{c}}}W] D^2 h\r)\r\|_{L^n\l(U\cap B_{\tilde{c}/2}\r)}\\
&\leq&C\l\{(1+\|u\|_{C^{\f{1}{2}}(\partial U\cap B_{c^2})})\|\Phi-(g^1)^{n-1}_{U\cap B_{\tilde{c}}}W\|^{\gz}_{L^n\l(U\cap B_{\tilde{c}}\r)}+\|f\|_{L^n\l(U\cap B_{c^2}\r)}\r\}.
\end{eqnarray*}
\end{lem}

With the above two lemmas, we next prove the basic power decay estimate for the distribution function of $D^2 u$, arguing as in the proof of \cite[Lemma 5.1]{NT}.

\begin{lem}\label{s3:lem NT5.1}
Assume $(\Omega,\phi, g, U)\in \mathcal{P}_{\lz,\,\Lz,\,\rz,\,\kappa,\az}$ with $0<\kappa<1$, $\phi\in C^1(\Omega)$ and $\mathrm{mosc}_U g\leq\ez$. Let $c$ be the constant in the definition of the class $\mathcal{P}$ and $\tilde{c}\leq \f{c^2}{2}$ be the constant given by Lemma \ref{s3:lem NT4.6}. Assume $u\in C(\Omega)\cap C^1(U)\cap W^{2,n}_{\mathrm{loc}}(U)$ is a solution of $\mathcal{L}_{\phi} u=f$ in $U$ such that
$$\|u\|_{L^\infty(U)}+\|u\|_{C^{2,\az}(\partial U\cap B_{c^2})}\leq 1$$
and
\begin{equation}\label{s3:NT5.11}
|u(x)|\leq C^*[1+d(x,x_0)^2]\quad\mathrm{in}\;\Omega\backslash U\;\mathrm{for\;some\;}x_0\in U\cap B_{\tilde{c}/4}.
\end{equation}
Then there exist $C=C(n,\lz,\Lz,\rz,\az,\kappa)>0$, $\tau=\tau(n,\lz,\Lz,\rz)>0$ and $N_0=N_0(n,\lz,\Lz,\rz,\az,\kappa,C^*)>0$ such that
\begin{eqnarray*}
|S_{\!\phi}(0,\l(\tilde{c}/2\r)^6)\cap G_{\!N}(u,\Omega)|\geq\l\{1-C\l(\f{\dz_0}{N}\r)^{\tau}-\ez^{\f{\tau'}{n}}\r\}|S_{\!\phi}(0,\l(\tilde{c}/2\r)^6)|,
\quad\forall N\geq N_0
\end{eqnarray*}
provided that $\|\Phi-(g^1)^{n-1}_{U\cap B_{\tilde{c}}}W\|_{L^n(U\cap B_{\tilde{c}})}\leq\l(\f{\tilde{c}}{2}\r)^{4}$. Here $\tau'=\tau'(n)$ is from Lemma \ref{s3:lem NT4.6}, $w$ is defined in $(vii)$ of the definition of $\mathcal{P}$ with $A:=[(g^1)_{U\cap B_{\tilde{c}}}]^{-1}$ and $W$ is the cofactor matrix of $D^2 w$, $\gz=\gz(n,\lz,\Lz,\rz,\az)$ is from Lemma \ref{s3:lem NT4.5}, and
$$\dz_0:=\l\{(1+\|u\|_{C^{\f{1}{2}}(\partial U\cap B_{c^2})})\|\Phi-(g^1)^{n-1}_{U\cap B_{\tilde{c}}}W\|^{\gz}_{L^n\l(U\cap B_{\tilde{c}}\r)}+\l(\fint_U |f|^n dx\r)^{\f{1}{n}}\r\}.$$
\end{lem}

Using Lemma \ref{s3:lem NT5.1}, the Localization Theorem, the interior power decay estimate in Lemmas \ref{s2:lem GN24.3},
\ref{s2:lem GN24.5} and the stablity of cofactor matrices at the boundary in Proposition \ref{s3:prop NT3.14},
we can prove power decay estimates when $\mathrm{det}\,D^2\phi$ satisfies a $\mathrm{VMO}$-type condition, following similar lines as in \cite[Lemma 5.2, Lemma 5.3]{NT}.

\begin{lem}\label{s3:lem NT5.2}
Given $0<\ez_0<1$ and $0<\az<1$. Assume $\Omega$ and $\phi$ satisfy \eqref{s3:g1}-\eqref{s3:g4}. Assume in addition that $\partial\Omega\in C^{2,\az}$ and $\phi\in C^{2,\az}(\partial\Omega)$. Let $u\in C^1(\Omega)\cap W^{2,n}_{\mathrm{loc}}(\Omega)$ be a solution of $\mathcal{L}_{\phi}u=f$ in $\Omega$ with $u=0$ on $\partial\Omega$ and $\|u\|_{L^\infty(\Omega)}\leq 1$. Then there exist $\ez=\ez(\ez_0,n,\lz,\Lz,\rz,\az)>0$, $c_1=c_1(n,\lz,\Lz,\rz,\az,\|\partial\Omega\|_{C^{2,\az}},\|\phi\|_{C^{2,\az}(\partial\Omega)})>0$ such that if
\begin{equation}\label{s3:NT5.21}
\sup_{x\in\overline{\Omega},\,t>0}\mathrm{mosc}_{S_{\!\phi}(x,t)}\,g\leq\ez,
\end{equation}
then for any $x\in\overline{\Omega}$ and $t\leq c_1$, we have
\begin{equation}\label{s3:NT5.22}
|G_{\!\f{N}{t}}(u,\Omega)\cap S_{\!\phi}(x,t)|\geq\l\{1-\ez_0-C\l(\f{\sqrt{t}}{N}\r)^{\tau}\|f\|_{L^n(\Omega)}^{\tau}\r\}|S_{\!\phi}(x,t)|,\quad\forall N\geq N_1,
\end{equation}
where $\tau=\tau(n,\lz,\Lz,\rz)$, $C,N_1>0$ depend only on $n,\lz,\Lz,\rz,\az$.
\end{lem}

\begin{lem}\label{s3:lem NT5.3}
Given $0<\ez_0<1$ and $0<\az<1$. Assume $\Omega$ and $\phi$ satisfy \eqref{s3:g1}-\eqref{s3:g4}. Assume in addition that $\Omega$ is uniformly convex, $\partial\Omega\in C^{2,\az}$ and $\phi\in C^{2,\az}(\partial\Omega)$. Let $u\in C^1(\Omega)\cap W^{2,n}_{\mathrm{loc}}(\Omega)$ be a solution of $\mathcal{L}_{\phi}u=f$ in $\Omega$ with $u=0$ on $\partial\Omega$. Then there exist $\ez=\ez(\ez_0,n,\lz,\Lz,\rz,\az)>0$, $c_2=c_2(n,\lz,\Lz,\rz,\az,\|\partial\Omega\|_{c^{2,\az}},\|\phi\|_{C^{2,\az}(\partial\Omega)})>0$ such that if
$$\sup_{x\in\overline{\Omega},\,t>0}\mathrm{mosc}_{S_{\!\phi}(x,t)}\,g\leq\ez,$$
then for any $x\in\overline{\Omega}$, $t\leq c_2$ and $S_{\!\phi}(x,t)\cap G_{\gz}(u,\Omega)\neq\emptyset$, we have
\begin{equation}\label{s3:NT5.32}
|G_{\!N\gz}(u,\Omega)\cap S_{\!\phi}(x,t)|\geq\l\{1-\ez_0-C(N\gz)^{-\tau}\l(\fint_{S_{\!\phi}(\tilde{x},\Theta t)}|f|^n dx\r)^{\f{\tau}{n}}\r\}|S_{\!\phi}(x,t)|,\quad\forall N\geq N_2,\tilde{x}\in S_{\!\phi}(x,t),
\end{equation}
where $\tau,\Theta$ depend only on $n,\lz,\Lz,\rz$; $C,N_2>0$ depend only on $n,\lz,\Lz,\rz,\az,\|\partial\Omega\|_{c^{2,\az}},\|\phi\|_{C^{2,\az}(\partial\Omega)}$ and the uniform convexity of $\Omega$.
\end{lem}

\section{Boundary estimate-Proof of the main theorems}\label{s7}

\noindent$\mathbf{Proof\;of\;Theorem\;2:}$
We assume that $\|u\|_{L^\infty(\Omega)}\leq 1$ and $\|f\|_{L^q(\Omega)}\leq\ez$, and we only need to prove that
\begin{equation}\label{s3:thm NT1.12}
\|D^2 u\|_{L^p(\Omega)}\leq C(n,\lz,\Lz,\Omega,p,q).
\end{equation}
Let $N_1,N_2$ be the large constants depending only on $n,\lz,\Lz,\Omega$ given by Lemmas \ref{s3:lem NT5.2} and \ref{s3:lem NT5.3} respectively. Denote $N_*:=\max\{N_1,N_2\}$. Let $c_1,c_2$ be the small constants depending only on $n,\lz,\Lz,\Omega$ given by Lemmas \ref{s3:lem NT5.2} and \ref{s3:lem NT5.3} respectively and denote $\hat{c}:=\min\{c_1,c_2\}$. Fix $M\geq N_*$ such that $1/M<\hat{c}$. Choose $0<\ez_0<1$ small such that
$$M^q\sqrt{2\ez_0}=\f{1}{2}.$$
Let $\ez=\ez(\ez_0,n,\lz,\Lz,\Omega)=\ez(n,\lz,\Lz,\Omega,q)$ be the smallest of the constants in Lemmas \ref{s3:lem NT5.2} and \ref{s3:lem NT5.3}. As in the proof of \cite[Theorem 1.1]{NT}, to obtain \eqref{s3:thm NT1.12}, we only need to estimate the term
$$\sum_{k=1}^\infty M^k\big|\Omega\backslash A^{\mathrm{loc}}_{\l(cM^{\f{k(q-p)}{2p}}\r)^{\f{-2}{n-1}}}\big|.$$
For this, we apply Theorem \ref{s3:thm NT3.5} with $\bz\rightsquigarrow M^k,\kappa\rightsquigarrow\f{q}{p}>1$ and find that
\begin{eqnarray*}
\sum_{k=1}^\infty M^k\big|\Omega\backslash A^{\mathrm{loc}}_{\l(cM^{\f{k(q-p)}{2p}}\r)^{\f{-2}{n-1}}}\big|
&\leq&C(\ez,n,\lz,\Lz,\Omega)\sum_{k=1}^\infty  M^{k\l(1+\f{(q-p)\ln\l(2\ez^{\dz_1/2}\r)}{2p(n-1)\ln C}\r)}<\infty
\end{eqnarray*}
if we choose $\ez=\ez(n,\lz,\Lz,\Omega,p,q)$ is small. Here $C,\dz_1>0$ are constants depending only on $n,\lz,\Lz,\Omega$.\\

\noindent$\mathbf{Proof\;of\;Theorem\;3:}$
Let $0<\ez_0=\ez_0(n,\lz,\Lz,p,q,\Omega)<1$ and $\ez=(\ez_0,n,\lz,\Lz,p,q,\Omega)>0$ be constants to be chosen later. We can assume that $\varphi=0, \|u\|_{L^\infty(\Omega)}\leq 1$ and $\|f\|_{L^q(\Omega)}\leq\ez$, and we only need to prove that
\begin{equation}\label{s3:thm NT1.2*}
\|D^2 u\|_{L^p(\Omega)}\leq C(n,\lz,\Lz,p,q,\Omega,g).
\end{equation}

Denote
\begin{equation*}
\eta_{g}(r,\Omega)=\sup_{\substack{x\in\overline{\Omega},\,t>0,\\\mathrm{diam}(S_{\!\phi}(x,t))\leq r}}\mathrm{mosc}_{S_{\!\phi}(x,t)}\,g.
\end{equation*}
Since $g\in\mathrm{VMO}(\Omega,\phi)$, then there exists $0<m<1$ depending only on $\ez$ and $g$ such that $\eta_{g}(r,\Omega)<\ez$ for any $0<r\leq m$. It follows that
\begin{equation}\label{s3:thm NT1.22}
\sup_{\substack{x\in\overline{\Omega},\,t>0,\\\mathrm{diam}(S_{\!\phi}(x,t))\leq m}}\mathrm{mosc}_{S_{\!\phi}(x,t)}\,g\leq\ez.
\end{equation}

Consider a point $y\in\partial\Omega$ and assume for simplicity that $y=0$. Assume that $\Omega$ satisfies \eqref{s3:1}, and $\phi(0)=0,\nabla\phi(0)=0$. By the Localization Theorem \ref{s3:thm NT2.2}, we have \begin{eqnarray}\label{s3:thm NT1.24}
\overline{\Omega}\cap B^+_{s^{2/3}}\subset S_{\!\phi}(0,s)\subset B_{s^{1/3}}
\end{eqnarray}
if $s\leq c_1=c_1(n,\lz,\Lz,\Omega)$. Hence if we choose $s$ such that $s\leq\min\l\{c_1,\l(\f{m}{2}\r)^3\r\}$, then \eqref{s3:thm NT1.22} implies that
\begin{equation}\label{s3:thm NT1.23}
\sup_{S_{\!\phi}(x,t)\subseteq S_{\!\phi}(0,s)}\mathrm{mosc}_{S_{\!\phi}(x,t)}\,g\leq\ez.
\end{equation}

Denote the rescaled functions and domains as in \eqref{s3:5} and \eqref{s3:6},
\begin{eqnarray*}
\phi_s(x)&=&\f{\phi(T_{\!s}^{-1}x)}{s}\quad\quad\mathrm{and}\quad\quad\Omega_s:=T_{\!s}(\Omega),\\
u_s(x)&=&u(T_{\!s}^{-1}x),\quad x\in\Omega_s,
\end{eqnarray*}
where $T_{\!s}:=s^{-\f{1}{2}}A_s$ and $A_s$ is the sliding given by Theorem \ref{s3:thm NT2.2}. Then
$$\mathrm{det}D^2\phi_s=g_s=(g^1_s)^n,\quad\quad 0<\lambda\leq g_s(x):=g(T_s^{-1}x)\leq\Lambda\quad\mathrm{in}\;\Omega_s,$$
$$\overline{\Omega_s}\cap B_k\subset U_s:=T_s(S_{\!\phi}(0,s))=S_{\!\phi_s}(0,1)\subset B_{k^{-1}}^+,$$
and by Proposition \ref{s3:prop NT2.12}, for any $0<\az<1$, we have
$$\l(\Omega_s,\phi_s, g_s, U_s\r)\in\mathcal{P}_{\lz,\,\Lz,\,\rz,\,C s^{1/2},\,\az}\subset\mathcal{P}_{\lz,\,\Lz,\,\rz,\,1/2,\,\az}$$
if $s\leq c_0$, where $c_0>0$ is a small constant depending only on $n,\lz,\Lz,\rz,\az,\|\partial\Omega\|_{C^{2,\az}}$ and $C=C(n,\lz,\Lz,\rz,\az)$.

Moreover, by \eqref{s3:thm NT1.23}, we have

\begin{eqnarray}\label{s3:thm NT1.26}
\sup_{S_{\!\phi_s}(x,t)\subset U_s\cap B_{k/2}}\mathrm{mosc}_{S_{\!\phi_s}(x,t)}\,g_s
&\leq&\sup_{S_{\!\phi}(T_s^{-1}x,\,ts)\subset S_{\!\phi}(0,s)}\mathrm{mosc}_{S_{\!\phi}(T_s^{-1}x,\,ts)}\,g\leq\ez.
\end{eqnarray}

Let $c=c(n,\lz,\Lz,\rz)$ be the small constant in $(vi)$ in the definition of $\l(\Omega_s,\phi_s, g_s, U_s\r)\in\mathcal{P}_{\lz,\,\Lz,\,\rz,\,1/2,\,\az}$. We claim that
\begin{eqnarray}\label{s3:thm NT1.27}
\|D^2 u_s\|_{L^p(S_{\!\phi_s}(0,c^{21}))}\leq C,
\end{eqnarray}
where $C>0$ depends only on $n,\lz,\Lz,p,q,\Omega$. Then back to $u$ and using a covering argument and the interior estimate Theorem \ref{s2:thm GN21.1}, 
the conclusion of Theorem 3 follows. 

To prove \eqref{s3:thm NT1.27} we establish the following two estimates which are the local versions of Lemmas \ref{s3:lem NT5.2} and \ref{s3:lem NT5.3} at the boundary respectively.

For any $x\in\overline{U_s}\cap B_{c^2}$, $t\leq c_1$ and $N\ge N_1$, we have
\begin{eqnarray}\label{s3:NT1.2a}
|S_{\!\phi_s}(x,t)\cap G_{\!\f{N}{t}}(u_s,U_s,\phi_s)|
&\geq&\l\{1-\ez_0-C\l(\f{\sqrt{t}}{N}\r)^\tau
\|f_s\|_{L^n(U_s)}^{\tau}\r\}|S_{\!\phi_s}(x,t)|.
\end{eqnarray}

For any $x\in\overline{U_s}\cap B_{c^6}$, $t\leq c_2$ and $N\ge N_2$, if $S_{\!\phi_s}(x,t)\cap G_{\!\gz}(u_s,U_s,\phi_s)\neq\emptyset$, then, for any $\tilde{x}\in S_{\!\phi_s}(x,t)$, we have
\begin{eqnarray}\label{s3:NT1.2b}
|S_{\!\phi_s}(x,t)\cap G_{\!N\gz}(u_s,U_s,\phi_s)|
&\geq&\l\{1-\ez_0-C\l(\f{1}{N\gz}\r)^{\tau}
\l(\fint_{S_{\!\phi_s}\l(\tilde{x},\Theta t\r)}|f_s|^n dx\r)^{\f{\tau}{n}}\r\}|S_{\!\phi_s}(x,t)|.
\end{eqnarray}
Here $c_1,c_2,N_1,N_2,\Theta$ are constants depending only on $n,\lz,\Lz,\Omega$.\\

The proof of these two results are similar to that of Lemmas \ref{s3:lem NT5.2}, \ref{s3:lem NT5.3}, but we use
localized versions of geometric properties of sections at the boundary (Dichotomy, Engulfing property, Volume estimates) instead of the
global ones. We omit the proof.

Using \eqref{s3:NT1.2a} and \eqref{s3:NT1.2b}, we are ready to prove \eqref{s3:thm NT1.27}. The proof of this is similar to
that of \cite[Theorem 1.1]{NT}, but we use the localized versions of covering theorem Lemma and
strong $p$-$p$ estimates for maximal functions instead of the global ones. We omit the proof.

LMAM, School of Mathematical Sciences,
 Peking University, Beijing, 100871,
 P. R. China

 Lin Tang,\quad
 E-mail address:  tanglin@math.pku.edu.cn

Qian Zhang,\quad
E-mail address: 1401110018@pku.edu.cn
\end{document}